\documentclass[oneside]{amsart}

\usepackage{geometry}

\usepackage{subfiles}
\usepackage{hyperref}
\usepackage[capitalise]{cleveref}
\usepackage{amsmath,amssymb}
\usepackage{tikz}
\usetikzlibrary{cd}
\usepackage{aligned-overset}

\newenvironment{spm}
  {\left( \begin{smallmatrix}}
  {\end{smallmatrix} \right) }

\DeclareMathOperator{\Bstd}{\textbf{B}^{std}}
\DeclareMathOperator{\Bstdv}{\textbf{B}^{std}_{\textbf{v}}}
\DeclareMathOperator{\Bcan}{\textbf{B}^{can}}
\DeclareMathOperator{\Bcanv}{\textbf{B}^{can}_{\textbf{v}}}
\DeclareMathOperator{\Bmon}{\textbf{B}^{mon}}
\DeclareMathOperator{\Bmonv}{\textbf{B}^{mon}_{\textbf{v}}}
\DeclareMathOperator{\KP}{KP}
\DeclareMathOperator{\Hom}{Hom}
\DeclareMathOperator{\GL}{GL}
\DeclareMathOperator{\Perv}{Perv}
\DeclareMathOperator{\IC}{IC}
\DeclareMathOperator{\pt}{pt}
\DeclareMathOperator{\Sbf}{\textbf{S}}
\DeclareMathOperator{\Irr}{Irr}
\DeclareMathOperator{\Ind}{Ind}

\DeclareMathOperator{\Semis}{Semis}
\DeclareMathOperator{\aff}{aff}
\DeclareMathOperator{\Diag}{Diag}
\DeclareMathOperator{\Id}{Id}
\DeclareMathOperator{\modu}{-mod}

\newcommand{\downdots}{\rotatebox[origin=c]{-26.9}{$\cdots$}}

\newtheorem{theorem}{Theorem}[section]
\newtheorem{corollary}[theorem]{Corollary}
\newtheorem{proposition}[theorem]{Proposition}
\newtheorem{lemma}[theorem]{Lemma}
\newtheorem{algorithm}[theorem]{Algorithm}
\newtheorem{remark}[theorem]{Remark}
\newtheorem{example}[theorem]{Example}

\newtheorem*{theorem*}{Theorem}
\newtheorem*{proposition*}{Proposition}

\date{\today}
\author{Jonas Antor}
\title{Canonical bases via pairing monomials}

\begin{document}

\begin{abstract}
For any quantum group of finite ADE type, we prove a new formula for the standard bilinear form evaluated at monomials. Combining this with ideas from the Lusztig-Shoji algorithm, we obtain a new algorithm that computes the canonical basis. In type A, the algorithm also computes composition multiplicities of standard modules for the affine Hecke algebra of $\GL_n$ and we explain how the algorithm can be extended to compute the dimensions of simple modules.
\end{abstract}

\maketitle
\section{Introduction}
\subsection*{Background}
Let $U^+$ be the positive part of a quantum group of finite ADE type. The algebra $U^+$ admits several standard bases of PBW type but it was observed by Lusztig \cite{lusztig1990canonical} that $U^+$ also admits a canonical basis. A similar situation arises when $U^+$ is replaced with the Hecke algebra which, apart from its standard basis, also admits another `canonical basis' also known as the Kazhdan-Lusztig basis \cite{kazhdan1979representations}.

More generally, canonical bases arise from geometric and algebraic categorifications where they correspond to classes of indecomposable objects in a Grothendieck group. These categorifications often allow to interpret the coefficients in the change of basis matrix between the standard and canonical basis as dimensions of stalks of perverse sheaves or composition multiplicities of simple modules in certain standard modules. The most famous example of this phenomenon is probably the Kazhdan-Lusztig conjecture \cite{kazhdan1979representations} proved in \cite{brylinski1981kazhdan, beilinson1981localisation} which realises composition multiplicities of Verma modules as certain coefficients in the Kazhdan-Lusztig basis of the corresponding Hecke algebra. Similarly, the canonical basis of $U^+$ encodes the stalks of simple perverse sheaves on a certain moduli space of quiver representations \cite{lusztig1990canonical} or alternatively composition multiplicities of certain modules over KLR algebras \cite{kato2014poincare,mcnamara2015finite,brundan2014homological}.

In the examples above, various algorithms and techniques have been developed to compute the canonical basis $\Bcan = \{ b_c \mid c \in C\}$ in terms of the standard basis $\Bstd = \{ a_c \mid c \in C \}$. The standard approach for this, which is already present in the work of Kazhdan and Lusztig, exploits that $b_c$ is the unique element which is invariant under a certain bar involution and satisfies
\begin{equation*}
b_c \in a_c + \sum_{c'< c} v^{-1}\mathbb{Z}[v^{-1}] a_{c'}
\end{equation*}
for an appropriate partial order $\le$ on $C$. The elements $b_c$ can then be computed by the following inductive procedure:
\begin{enumerate}
\item Find bar invariant elements $\tilde{b}_c$ ($c \in C$) such that $\tilde{b}_c \in  a_c + \sum_{c'<c} \mathbb{Z}[v,v^{-1}] a_{c'}$;
\item Subtract $\mathbb{Z}[v+v^{-1}]$-multiples of the $b_{c'}$ ($c'< c$) from $\tilde{b}_c$ to ensure that the coefficients of the $a_{c'}$ lie in $v^{-1} \mathbb{Z}[v^{-1}]$. The resulting element is $b_c$.
\end{enumerate}
The main ingredient needed to make this algorithm work is an efficient way to express the elements $\tilde{b}_c$ in terms of the standard basis. In the case of the Hecke algebra, this can be achieved using explicit knowledge of the multiplication rules on the standard basis. A similar approach was carried out in the quantum group setting in \cite{de2002constructing} which also relies on explicit computations with products of PBW basis elements.

\subsection*{The standard bilinear form}
In this paper, we develop a new algorithm to compute the canonical basis of $U^+$ which does not require explicit knowledge of the multiplication. The algorithm is inspired by the Lusztig-Shoji algorithm \cite{shoji1983green,lusztig1986character} which computes stalks of perverse sheaves on the nilpotent cone. At the heart of the Lusztig-Shoji algorithm lies a Gram-Schmidt type procedure which is applied to a certain Ext-pairing on perverse sheaves (see \cite{achar2011green} for more details on this perspective). In the quantum groups setting, the replacement for the Ext-pairing is the standard bilinear form
\begin{equation*}
(\cdot, \cdot) : U^+ \times U^+ \rightarrow \mathbb{Q}(v)
\end{equation*}
defined in \cite[1.2.3]{lusztig2010introduction}. One of our main results is a new explicit formula for this bilinear form. The formula is combinatorial but its proof relies on Lusztig's geometric categorification of $U^+$. In principle, the bilinear form on $U^+$ can also be computed from its recursive definition (see also the closed formula in \cite{reineke2001monomials}) but our new formula has the advantage that it applies to more general geometric situations (\cref{theorem: Ext formula}). As an application of this, we will also develop an algorithm that computes the dimensions of the simple modules for the affine Hecke algebra of $\GL_n$.

To state our new formula, we need to briefly recall the construction of monomials in $U^+$. The algebra $U^+$ has standard generators $E_i$ ($i \in I$) and admits a grading
\begin{equation*}
U^+ = \bigoplus_{\textbf{v} \in \mathbb{N}[I]} U^+_{\textbf{v}}
\end{equation*}
where $E_i$ has degree $i$. For the purpose of computing the canonical basis, it is sufficient to restrict our attention to a fixed homogeneous component $U^+_{\textbf{v}}$. For any pair of sequences $\textbf{i} = (i_1, ..., i_m) \in I^m$, $\textbf{a} = (a_1, ..., a_m) \in \mathbb{N}^m$ of weight $\textbf{v}$ (i.e. $\sum_{i_k = j} a_k = \textbf{v}_j$ for all $j \in I$), one can consider the monomial
\begin{equation*}
E_{\textbf{i}, \textbf{a}} := E_{i_1}^{(a_1)} \cdot ... \cdot E_{i_m}^{(a_m)} \in U^+_{\textbf{v}}.
\end{equation*}
Choose an orientation of the Dynkin quiver associated to $U^+$ and let
\begin{equation*}
R_{\textbf{v}} = \bigoplus_{i \rightarrow j} \Hom(\mathbb{C}^{\textbf{v}_i}, \mathbb{C}^{\textbf{v}_j})
\end{equation*}
be the space of $\textbf{v}$-dimensional representations of this quiver. Let $\GL_{\textbf{v}}$ be the corresponding change of basis group with Weyl group $S_{\textbf{v}}$ and let $B \subset \GL_{\textbf{v}}$ be a Borel subgroup. We can associate to the pair $(\textbf{i}, \textbf{a})$ a certain subspace
\begin{equation*}
R_{\textbf{i}, \textbf{a}} \subset R_{\textbf{v}}
\end{equation*}
consisting of those representations that fix the unique $B$-stable flag of $\mathbb{C}^{\textbf{v}}$ of type $(\textbf{i}, \textbf{a})$ (see \cref{section: geometric categorification} for more details). We set $[\textbf{a}]! := [a_1]! \cdot ... \cdot [a_m]!$ where $[a_i]! \in \mathbb{Q}(v)$ is the quantum factorial. The following theorem is our new formula for the bilinear form on $U^+$.
\begin{theorem*}[\cref{cor: geometric counting formula for inner product}]
Let $(\textbf{i}, \textbf{a})$ and $(\textbf{i}', \textbf{a}')$ be of weight $\textbf{v}$. Then
\begin{equation*}
(E_{\textbf{i}, \textbf{a}}, E_{\textbf{i}', \textbf{a}'}) = v^{- \dim R_{\textbf{i}, \textbf{a}} - \dim R_{ \textbf{i}', \textbf{a}'}}  \frac{\prod_{i \in I} (1-v^{-2})^{-\textbf{v}_i}}{ [\textbf{a}]! [\textbf{a}']!} \sum_{w \in S_{\textbf{v}}} v^{2(l(w) + \dim R_{\textbf{i}, \textbf{a}} \cap  {}^w R_{\textbf{i}', \textbf{a}'})}.
\end{equation*}
\end{theorem*}
The pairing $(E_{\textbf{i}, \textbf{a}}, E_{\textbf{i}', \textbf{a}'})$ can also be interpreted as the graded dimension of the $\Hom$-space between
\begin{itemize}
\item certain semisimple complexes of perverse sheaves on $R_{\textbf{v}}$ \cite{lusztig2010introduction} (see also \cref{thm: geometric categorification of qua group}),
\item certain projective modules over a KLR algebra \cite[3.4]{khovanov2009diagrammatic}.
\end{itemize}
Thus, the formula above also computes these graded dimensions.

\subsection*{The algorithm}
We now explain our new algorithm. Let $\Bstdv$ (resp. $\Bcanv$) be the standard (resp. canonical) basis of $U^+_{\textbf{v}}$. Pick a basis $\Bmonv$ of $U^+_{\textbf{v}}$ consisting of monomials $E_{\textbf{i}, \textbf{a}}$ such that the change of basis matrix between $\Bcanv$ and $\Bmonv$ is lower unitriangular. This is always possible by \cite{lusztig1990canonical,reineke2001feigin,reineke2003quivers} (see also the discussion after \cref{lem: existence of res of singularities}). Let $\Psi$ be the matrix representing $(\cdot, \cdot)$ with respect to the monomial basis $\Bmonv$. The following proposition shows that the canonical basis is uniquely determined by $\Psi$.
\begin{proposition*}(\cref{section: linear algebra part})
There are unique matrix decompositions
\begin{equation*}
\Psi = L D L^T, \quad  L =QP
\end{equation*}
such that
\begin{enumerate}
\item $L, P$ and $Q$ are lower unitriangular;
\item $D$ is diagonal;
\item $P$ has entries in $v^{-1}\mathbb{Z}[v^{-1}]$ below the diagonal;
\item $Q$ has entries in $\mathbb{Z}[v+v^{-1}]$.
\end{enumerate}
The matrices $L, P, Q$ and $\Lambda$ also have the following alternative interpretations:
\begin{itemize}
\item $P$ is the change of basis matrix between $\Bstdv$ and $\Bcanv$;
\item $Q$ is the change of basis matrix between $\Bcanv$ and $\Bmonv$;
\item $L$ is the change of basis matrix between $\Bstdv$ and $\Bmonv$;
\item $D$ is the matrix representing $(\cdot, \cdot)$ with respect to $\Bstdv$.
\end{itemize}
\end{proposition*}
In linear algebra, the decomposition $\Psi = LD L^T$ is often referred to as the LDLT decomposition. It can be computed using the Gram-Schmidt orthogonalization algorithm. The decomposition $L=QP$ can be computed by a basic row elimination algorithm. Moreover, $\Psi$ can be computed using our new formula for the standard bilinear form on $U^+_{\textbf{v}}$. Combining these results, we arrive at the following algorithm to compute the canonical basis (see \cref{section: examples} for examples that apply this algorithm).
\begin{algorithm}\label{algorithm: new algo}Input: $\textbf{v}$ (weight vector), Output: $P,Q$ (change of basis matrix between $\Bstdv$ and $\Bcanv$ resp. $\Bcanv$ and $\Bmonv$).
\begin{itemize}
\item[Step 1:] Compute the matrix $\Psi$ using the formula from \cref{cor: geometric counting formula for inner product};
\item[Step 2:] Compute the decomposition $\Psi = LDL^T$ using the Gram-Schmidt procedure;
\item[Step 3:] Compute the decomposition $L = QP$ using row elimination.
\end{itemize}
\end{algorithm}
This algorithm can also be viewed through the lens of the classical Kazhdan-Lusztig algorithm: The monomial basis is an approximation of the canonical basis. Step 2 in the algorithm above expresses the monomial basis in terms of the standard basis. Computing the decomposition $L=QP$ in step 3 is then just a variation of the elimination procedure in the Kazhdan-Lusztig algorithm.
\subsection*{Affine Hecke algebras}
We now want to explain how the algorithm above can be extended to study the representation theory of the affine Hecke algebra $\mathcal{H}^{\aff}$ associated to $\GL_n$. The connection with the quantum groups situation is as follows: The representation theory of $\mathcal{H}^{\aff}$ can be studied using perverse sheaves on certain fixed-point varieties of the nilpotent cone $\mathcal{N}^{(s,q_0)}$ \cite{kazhdan1987proof,chriss2009representation}. For example, there are certain standard modules of $\mathcal{H}^{\aff}$ whose composition multiplicities can be interpreted as dimensions of stalks of perverse sheaves on $\mathcal{N}^{(s,q_0)}$. It turns out that when $q_0$ is not a root of unity, the space $\mathcal{N}^{(s,q_0)}$ is isomorphic to a moduli space $R_{\textbf{v}}$ of type $A$ quiver representations (or products thereof). Via Lusztig's geometric categorification of $U^+$, the stalks of perverse sheaves on $\mathcal{N}^{(s,q_0)} \cong R_{\textbf{v}}$ can thus be interpreted as the coefficients of the change of basis matrix between the bases $\Bstdv$ and $\Bcanv$ of $U^+_{\textbf{v}}$. Hence, \cref{algorithm: new algo} also computes these stalks and thus also the composition multiplicities for the standard modules.

Another quantity that admits a geometric interpretation in this setting are the dimensions of the simple $\mathcal{H}^{\aff}$-modules. It can be shown that there are certain semisimple complexes $\Sbf^{(s,q_0)} \in D^b_c(\mathcal{N}^{(s,q_0)})$ such that the dimensions of simple $\mathcal{H}^{\aff}$-modules are given by the multiplicities of the simple perverse sheaves in $\Sbf^{(s,q_0)}$ \cite{chriss2009representation}. Under the isomorphism $\mathcal{N}^{(s,q_0)} \cong R_{\textbf{v}}$ and Lusztig's geometric categorification of $U^+$, the sheaves $\Sbf^{(s,q_0)}$ correspond to certain elements in $U^+_{\textbf{v}}$. Decomposing the sheaves $\Sbf^{(s,q_0)}$ into simples can then be reduced to a certain Ext-computation which we can perform using a variation of our formula for the standard bilinear form. We spell out this extended algorithm for computing dimensions of simple modules in more detail in \cref{algorithm: dimensions of simple modules}. As an example, we compute the dimensions of certain simple modules for the affine Hecke algebra of $\GL_4$ in \cref{example: dimensions of simple for v equals 121}.

\subsection*{Acknowledgements}
The author would like to thank Kevin McGerty and Dan Ciubotaru for many helpful discussions. This research was financially supported by the Lincoln-Kingsgate Graduate Scholarship.

\section{The linear algebraic skeleton}\label{section: linear algebra part}
Let $V$ be a finite-dimensional $\mathbb{Q}(v)$-vector space equipped with the following data:
\begin{itemize}
\item a bar involution $\overline{\cdot} : V \rightarrow V$ which is $\mathbb{Q}(v)$-antilinear;
\item a standard basis $\Bstd= \{a_c \mid  c \in C \}$ indexed by a partially ordered set $(C, \le) $ satisfying
\begin{equation*}
\overline{a}_c \in a_c + \sum_{c' < c} \mathbb{Z}[v,v^{-1}]\cdot a_{c'}.
\end{equation*}
\end{itemize}
Then by \cite[24.2.1]{lusztig2010introduction} there is a unique \textit{canonical basis}
\begin{equation*}
\Bcan = \{ b_c \mid c \in C\}
\end{equation*}
of $V$ satisfying
\begin{itemize}
\item $\overline{b}_c = b_c$ for all $c \in C$;
\item $b_c = a_c + \sum_{c' < c} p_{cc'}(v) a_{c'}$ for some $p_{cc'}(v) \in v^{-1}\mathbb{Z}[v^{-1}]$.
\end{itemize}
If we refine the partial order on $C$ to a total order, the last condition implies that the change of basis matrix
\begin{equation*}
P := (p_{cc'}(v))_{c,c'\in C}
\end{equation*}
between $\Bstd$ and $\Bcan$ is lower unitriangular with entries in $v^{-1}\mathbb{Z}[v^{-1}]$ below the diagonal. We want to compute that canonical basis in terms of the standard basis or equivalently the matrix $P$. For this, we assume that we are given the following additional data: Let
\begin{equation*}
( \cdot, \cdot ) : V \times V \rightarrow \mathbb{Q}(v)
\end{equation*}
be a non-degenerate symmetric bilinear form such that $\Bstd$ is orthogonal with respect to $(\cdot, \cdot)$. Consider the matrices
\begin{align*}
D &:= ((a_c, a_{c'}))_{c,c' \in C} \\
\Omega &:= ((b_c, b_{c'}))_{c,c' \in C} 
\end{align*}
representing $(\cdot, \cdot)$ with respect to $\Bstd$ and $\Bcan$. By change of basis, we have
\begin{equation}\label{eq: LDLT decomposition of Omega}
\Omega = P D P^T
\end{equation}
where $\Omega$ is symmetric and invertible, $P$ is lower unitriangular and $D$ is diagonal. In linear algebra, this is called the LDLT decomposition of $\Omega$. The LDLT decomposition of any symmetric invertible matrix is unique (if it exists) and it can be computed using the Gram-Schmidt orthogonalization algorithm.
\begin{corollary}
We can compute $\Bcan$ from the matrix $\Omega$, i.e. from the values of the bilinear form on $\Bcan$.
\end{corollary}
In applications, computing the bilinear form on $\Bcan$ is a non-trivial task. However, there are sometimes certain `monomial elements' in $V$ for which the bilinear form can be computed explicitly. Assume that we are given the following additional data: Let
\begin{equation*}
\Bmon = \{ m_c \mid c \in C \}
\end{equation*}
be another ordered basis called the 'monomial basis' satisfying
\begin{itemize}
\item $\overline{m}_c = m_c$ for all $c \in C$;
\item $m_c = b_c + \sum_{c' < c} q_{cc'} b_{c'}$ for some $q_{cc'} \in \mathbb{Z}[v,v^{-1}]$
\end{itemize}
The last condition implies that the change of basis matrix
\begin{equation*}
Q := (q_{cc'})_{c,c' \in C}
\end{equation*}
between $\Bcan$ and $\Bmon$ is lower unitriangular with entries in $\mathbb{Z}[v,v^{-1}]$. This is equivalent to the change of basis matrix
\begin{equation*}
L := QP
\end{equation*}
between $\Bstd$ and $\Bmon$ being lower unitriangular with entries in $\mathbb{Z}[v,v^{-1}]$. Since the elements of $\Bcan$ and $\Bmon$ are bar invariant, we have
\begin{equation}\label{eq: matrix Q is bar invariant}
\overline{q}_{cc'} = q_{cc'}
\end{equation}
for all $c, c' \in C$. Thus the entries of $Q$ actually lie in $\mathbb{Z}[v+v^{-1}]$. Let
\begin{equation*}
\Psi:= ((m_c, m_{c'}))_{c,c' \in C}
\end{equation*}
be the matrix representing the bilinear form $(\cdot, \cdot)$ with respect to the basis $\Bmon$. Then by change of basis we get
\begin{equation}\label{eq: LDLT decomposition of Psi}
\Psi = L D L^T
\end{equation}
which is the LDLT decomposition of $\Psi$. Hence, $L$ and $D$ can be computed from $\Psi$. Note that we also have
\begin{equation}\label{eq: Psi is Q omega Q t}
\Psi = Q \Omega Q^T.
\end{equation}
The following lemma says that $P$ and $Q$ are uniquely determined by $L$ (and hence also by $\Psi$).
\begin{lemma}
The matrices $P$ and $Q$ are the unique matrices satisfying
\begin{itemize}
\item $L = QP$;
\item $Q$ is lower unitriangular with entries in $\mathbb{Z}[v+v^{-1}]$;
\item $P$ is lower unitriangular with entries in $v^{-1} \mathbb{Z}[v^{-1}]$ below the diagonal;
\end{itemize}
Moreover, $P$ and $Q$ can be computed from $L$ by a row elimination algorithm.
\end{lemma}
\begin{proof}
Assume we have $L = Q_1P_1 = Q_2P_2$ satisfying the conditions above. Then $P_1P_2^{-1} = Q_2^{-1}Q_1$ is a lower triangular matrix whose entries below the diagonal lie in $\mathbb{Z}[v+v^{-1}] \cap v^{-1}\mathbb{Z}[v^{-1}] = \{ 0 \}$. Hence, $P_1 = P_2$ and $Q_1 = Q_2$ which proves uniqueness. To compute $P$ and $Q$, note that we can transform $L$ into a matrix with entries in $v^{-1}\mathbb{Z}[v^{-1}]$ below the diagonal using only elementary row operations with coefficients in $\mathbb{Z}[v+v^{-1}]$. The result of this process is the matrix $P$ and the product of all the row operation is the matrix $Q^{-1}$.
\end{proof}
We can summarize the results of this section as follows.
\begin{corollary}\label{cor: can base can be computed from psi}
The matrix $P$ can be computed from the matrix $\Psi$ by computing the unique matrix decompositions $\Psi = LDL^T$ and $L = QP$.
\end{corollary}
Hence, the problem of computing the canonical basis is reduced to finding a nice monomial basis for which the matrix $\Psi$ can be computed explicitly.
\section{Quantum groups}\label{section: quantum groups}
In this section we recall a few well-known facts about the positive half of quantum groups and its standard, canonical and monomial bases. We denote by $\mathbb{N}$ the set of all natural numbers including zero. For any $k \in \mathbb{N}$, we write $[k]! \in \mathbb{Q}(v)$ for the quantum factorial.
\subsection{PBW bases}
Let $\Phi$ be a simply-laced root system with Weyl group $W$. Choose a set of positive roots $\Phi^+ \subset \Phi$ and denote the corresponding simple roots by $\textbf{e}_i$ ($i \in I$). The positive half of the quantum group associated to $\Phi$ is the $\mathbb{Q}(v)$-algebra $U^+ = U^+(\Phi)$ with generators $E_i$ ($i \in I$) and relations
\begin{align*}
E_iE_j - E_j E_i & = 0 \quad \text{ if  } \langle \textbf{e}_i, \check{\textbf{e}}_j \rangle = 0, \\
E_i^2 E_j -  (v+v^{-1})E_i E_j E_i + E_j E_i^2 & = 0 \quad \text{ if } \langle \textbf{e}_i, \check{\textbf{e}}_j \rangle = -1.
\end{align*}
This algebra admits a grading
\begin{equation*}
U^+ = \bigoplus_{\textbf{v} \in \mathbb{N}[I]} U^+_{\textbf{v}}
\end{equation*}
where $E_i$ has degree $i$. For any $x \in U^+$ we define
\begin{equation*}
x^{(k)} := \tfrac{x^k}{[k]!}.
\end{equation*}
The bar involution is the ring homomorphism
\begin{equation*}
\overline{\cdot} : U^+ \rightarrow U^+
\end{equation*}
with $\overline{v} = v^{-1}$ and $\overline{E}_i = E_i$. This restricts to an involution
\begin{equation*}
\overline{\cdot} : U^+_{\textbf{v}} \rightarrow U^+_{\textbf{v}}
\end{equation*}
on each homogeneous component. Fix a reduced expression 
\begin{equation*}
w_0 =s_{i_1} s_{i_2} ... s_{i_N}
\end{equation*}
of the longest element of $W$. This induces a total order on the set of positive roots
\begin{equation*}
\Phi^+ = \{ \gamma_1 < \gamma_2 < ... < \gamma_N\}
\end{equation*}
where
\begin{equation*}
\gamma_j := s_{i_1} s_{i_2} ... s_{i_{j-1}} ( \textbf{e}_{i_j}).
\end{equation*}
Moreover, using braid group actions, one can define for each $\gamma \in \Phi^+$ a distinguished element
\begin{equation*}
E_{\gamma} \in U^+
\end{equation*}
which also depends on the choice of reduced expression for $w_0$ \cite[\S 40]{lusztig2010introduction}. For any $\textbf{v} \in \mathbb{N}[I]$ let
\begin{equation*}
\KP(\textbf{v}) = \{ \textbf{c} : \Phi^+ \rightarrow \mathbb{N} \mid \sum_{\gamma \in \Phi^+} \textbf{c}(\gamma)\cdot \gamma = \sum_{i \in I} \textbf{v}_i \textbf{e}_i \}
\end{equation*}
be the set of Kostant partitions of $\textbf{v}$. For each $\textbf{c} \in \KP(\textbf{v})$ the corresponding PBW element is defined as
\begin{equation*}
E_{\textbf{c}} := E_{\gamma_1}^{(c(\gamma_1))} E_{\gamma_2}^{(c(\gamma_2))}  ... E_{\gamma_N}^{(c(\gamma_N))} \in U^+_{\textbf{v}}.
\end{equation*}
These elements form a basis
\begin{equation*}
\Bstdv  = \{ E_{\textbf{c}} \mid \textbf{c} \in \KP(\textbf{v}) \}
\end{equation*}
of $U^+_{\textbf{v}}$ known as the standard basis (or PBW basis) \cite[40.2.2]{lusztig2010introduction}. There is a coproduct $U^+ \rightarrow U^+ \otimes U^+$ which equips $U^+$ with the structure of a (twisted) bialgebra. Moreover there is a non-degenerate symmetric bilinear form
\begin{equation}\label{eq: bilinear form on qua groups}
( \cdot, \cdot) : U^+ \times U^+ \rightarrow U^+
\end{equation}
defined in \cite[1.2.3]{lusztig2010introduction} which is given on generators by
\begin{equation*}
( E_i, E_j ) = \delta_{i,j} (1-v^{-2})^{-1}
\end{equation*}
and adjoins the product in $U^+$ with the coproduct. This bilinear form satisfies
\begin{equation*}
( U^+_{\textbf{v}}, U^+_{\textbf{w}}) = 0
\end{equation*}
for any $\textbf{v} \neq \textbf{w}$ and thus the restricted bilinear form
\begin{equation*}
(\cdot, \cdot): U^+_{\textbf{v}} \times U^+_{\textbf{v}} \rightarrow \mathbb{Q}(v)
\end{equation*}
is still non-degenerate.
\begin{lemma}\cite[38.2.3]{lusztig2010introduction} The standard basis $\Bstdv$ on $U^+_{\textbf{v}}$ is orthogonal with respect to this bilinear form.
\end{lemma}
\subsection{Canonical and monomial bases}
Fix an orientation of the Dynkin diagram of $\Phi$, i.e. a quiver $Q$ whose underlying graph is the Dynkin diagram. The structure of $U^+$ is closely related to the quiver $Q$. For any $\textbf{v} \in \mathbb{N}[I]$ we denote by
\begin{equation*}
R_{\textbf{v}} = R_{\textbf{v}}(Q) = \bigoplus_{i \rightarrow j} \Hom(\mathbb{C}^{\textbf{v}_i}, \mathbb{C}^{\textbf{v}_j} )
\end{equation*}
the space of $\textbf{v}$-dimensional representations of 
$Q$. This comes with a canonical change of basis action of
\begin{equation*}
\GL _{\textbf{v}} = \prod_{i \in I} \GL_{\textbf{v}_i}.
\end{equation*} 
By Gabriel's theorem there is a canonical bijection
\begin{equation}\label{eq: kostant partitions parametrize orbits}
\begin{aligned}
\KP(\textbf{v}) & \overset{1:1}{\leftrightarrow} \{ \GL _{\textbf{v}} \text{-orbits on } R_{\textbf{v}} \} \\
 \textbf{c} & \mapsto \mathcal{O}_{\textbf{c}}.
\end{aligned}
\end{equation}
This induces a partial order $\le$ on $\KP(\textbf{v})$ where
\begin{equation*}
\textbf{c} < \textbf{c}' \quad  \text{ if } \quad \dim \mathcal{O}_{\textbf{c}} < \dim \mathcal{O}_{\textbf{c}'} .
\end{equation*}
From now on, we assume that the reduced expression $w_0 = s_{s_1} s_{i_2} ... s_{i_N}$ is adapted to $Q$ in the sense of \cite[4.7]{lusztig1990canonical} (such a reduced expression always exists by \cite[4.12]{lusztig1990canonical}). Then by \cite[7.9]{lusztig1990canonical} we have
\begin{equation}\label{eq: lower unitriangularity of bar involution}
\overline{E}_{\textbf{c}} \in E_{\textbf{c}} + \sum_{\textbf{c}' < \textbf{c}} \mathbb{Z}[v,v^{-1}] E_{\textbf{c}'}
\end{equation}
for each $\textbf{c} \in \KP(\textbf{v})$. Hence, by \cite[24.2.1]{lusztig2010introduction} there is a canonical basis of $U^+_{\textbf{v}}$
\begin{equation*}
\Bcanv =  \{ b_{\textbf{c}} \mid \textbf{c} \in \KP(\textbf{v}) \}.
\end{equation*}
It turns out that this basis is independent of the choice reduced expression and orientation \cite[2.3, 3.2]{lusztig1990canonical} (however the partial order on $\KP(\textbf{v})$ does depend on the choice of orientation).
\begin{remark}
We can also replace the orbit dimension partial order $\le$ in \eqref{eq: lower unitriangularity of bar involution} with the orbit closure partial order $\le_{cl}$ where
\begin{equation*}
\textbf{c} \le_{cl} \textbf{c}' \quad \text{ if }  \quad \mathcal{O}_{\textbf{c}} \subset \overline{\mathcal{O}}_{\textbf{c}'}.
\end{equation*}
We work with $\le$ since it is easier to compute in examples.
\end{remark}
Next, we discuss the monomial basis. Given a non-negative integer $m$ and a pair of sequences
\begin{align*}
\textbf{i} &= (i_1, ..., i_m) \in I^m \\
\textbf{a} &= (a_1, ..., a_m) \in \mathbb{N}^m
\end{align*}
we say that $(\textbf{i}, \textbf{a})$ has weight $\textbf{v}$ if
\begin{equation*}
\textbf{v}_j = \sum_{i_k = j} a_k
\end{equation*}
for all $j \in I$. For any pair $(\textbf{i}, \textbf{a})$ of weight $\textbf{v}$ we can define a corresponding monomial
\begin{equation*}
E_{\textbf{i}, \textbf{a}} := E_{i_1}^{(a_1)} E_{i_2}^{(a_2)} ... E_{i_m}^{(a_m)} \in U^+_{\textbf{v}}.
\end{equation*}
\begin{lemma}\label{lem: existence of monomial basis}\cite[7.8]{lusztig1990canonical}, \cite[4.2]{reineke2001feigin}
We can find $m \in \mathbb{N}$, $\textbf{i} \in I^m$, and for each $\textbf{c} \in \KP(\textbf{v})$ a sequence $\textbf{a}_{\textbf{c}} \in \mathbb{N}^m$ such that $(\textbf{i}, \textbf{a}_{\textbf{c}})$ has weight $\textbf{v}$ and
\begin{equation*}
E_{\textbf{i}, \textbf{a}_c} \in E_{\textbf{c}} + \sum_{\textbf{c}'< \textbf{c}} \mathbb{Z}[q,q^{-1}]  E_{\textbf{c}'}
\end{equation*}
for all $\textbf{c} \in \KP(\textbf{v})$.
\end{lemma}
Thus, if we pick $\textbf{i}$ and for each $\textbf{c} \in \KP (\textbf{v})$ a sequence $\textbf{a}_{\textbf{c}}$ as above, we obtain a monomial basis
\begin{equation*}
\Bmonv = \{ E_{\textbf{i}, \textbf{a}_{\textbf{c}}} \mid \textbf{c} \in \KP(\textbf{v}) \}
\end{equation*}
of $U^+_{\textbf{v}}$. By \cref{cor: can base can be computed from psi} the canonical basis can then be computed from the matrix
\begin{equation}\label{eq: def of psi in qua group case}
\Psi = ((E_{\textbf{i}, \textbf{a}_{\textbf{c}}}, E_{\textbf{i}, \textbf{a}_{\textbf{c}'}}))_{\textbf{c}, \textbf{c}' \in \KP(\textbf{v})}.
\end{equation}
In principle, the entries of the matrix $\Psi$ can be computed from the recursive formula in \cite[1.2.3]{lusztig2010introduction}. We will provide an alternative formula for $\Psi$ in terms of quiver representations in the following section.
\begin{remark}
There is an explicit procedure \cite[4.2]{reineke2001feigin} that constructs the sequences $\textbf{a}_c$, depending on the choice of a 'directed partition' on the set of indecomposable representations of $Q$. We will give explicit examples in type A in \cref{section: examples}.
\end{remark}
\section{A geometric formula for the bilinear pairing}
In this section, we derive a new formula for the standard bilinear form of $U^+$ evaluated at monomials. The proof relies on geometric techniques.
\subsection{Notation}
Let $G$ be a connected complex algebraic group and $X$ a $G$-variety. We denote by $H^*_G(X)$ the equivariant cohomology of $X$ with complex coefficients and by $H_*^G(X)$ the equivariant Borel-Moore homology. For $G=\{e\}$ this recovers ordinary cohomology $H^*(X)$ and Borel-Moore homology $H_*(X)$. If $X$ is smooth and connected there is the Poincar\'e duality isomorphism
\begin{equation*}
H_k^G(X) \overset{PD}{\cong} H_G^{2\dim X -k}(X).
\end{equation*}
We will make use of the following standard result about Borel-Moore homology.
\begin{lemma}\label{parity vanishing direct sum of Borel Moore}
Let $X = \bigsqcup_{j \in J} X_j$ be a finite partition into locally closed $G$-stable subsets such that $H_k^G(X_j)$ vanishes for $k$ odd for all $j \in J$. Assume that we can pick a total order on $J$ such that $X_{\le j} = \bigsqcup_{j' \le j} X_{j'}$ is closed in $X$ for all $j \in J$. Then $H^G_k(X) \cong \bigoplus_{j \in J} H^G_k(X_j)$ for all $k$.
\end{lemma}
\begin{proof}
This follows by a standard induction argument using the long exact sequence of Borel-Moore homology (see the proof of \cite[B.3-Lemma 6]{fulton1997young} for more details on this type of argument).
\end{proof}
Denote by
\begin{equation}\label{eq: def of poincare poly}
\begin{aligned}
\chi_v^G(X) & := v^{\dim X} \sum_{k \ge 0} \dim H_G^k(X)v^{-k}  \in \mathbb{Z}((v^{-1}))\\
\chi_v(X) & := v^{\dim X} \sum_{k \ge 0} \dim H^k(X)v^{-k} \in \mathbb{Z}[v,v^{-1}]
\end{aligned}
\end{equation}
the (normalized) Poincar\'e polynomials of $X$. We write $D^b_G(X)$ for the equivariant constructible derived category with complex coefficients and $\textbf{1}_X \in D^b_G(X)$ for the constant sheaf. All sheaf functors will be understood to be derived. The Verdier duality functor will be denoted by $\mathbb{D}$. For any $\mathcal{F}, \mathcal{G} \in D^b_G(X)$ we define
\begin{equation*}
\Hom^k(\mathcal{F}, \mathcal{G}) := \Hom_{D^b_G(X)}(\mathcal{F}, \mathcal{G}[k]).
\end{equation*}
We write $\Perv_G(X) \subset D^b_G(X)$ for the category of equivariant perverse sheaves and $\Irr(\Perv_G(X))$ for the set of simple equivariant perverse sheaves. If $X$ is smooth with connected components $X_1, ..., X_n$, we write
\begin{equation*}
\mathcal{C}_X := \textbf{1}_{X_1}[\dim X_1] + ...+ \textbf{1}_{X_n} [\dim X_n] \in \Perv_G(X)
\end{equation*}
for the constant perverse sheaf. Denote by $\Semis_G(X) \subset D^b_G(X)$ the additive category of semisimple complexes, i.e. the full subcategory whose objects are direct sums of shifts is simple perverse sheaves. Then the split Grothendieck group
\begin{equation*}
K_{\oplus}(X) := K_{\oplus}(\Semis_G(X))
\end{equation*}
is a free $\mathbb{Z}[v,v^{-1}]$-module (where $v$ acts by the homological shift $[1]$) with a basis given by the classes of simple equivariant perverse sheaves. If $X$ is smooth and connected with structure map $a: X \rightarrow \pt$ we have
\begin{equation}\label{eq: Groth group interpretation of poincare poly}
[a_* \mathcal{C}_X] = \chi_v(X) \cdot [\mathcal{C}_{\pt}]
\end{equation}
in $K_{\oplus}(\pt)$. For any Laurent polynomial
\begin{equation*}
f(v) = \sum_{k \in \mathbb{Z}} a_k v^k \in \mathbb{N}[v,v^{-1}]
\end{equation*}
and $\mathcal{F} \in \Semis_G(X)$ we define
\begin{equation*}
f(v) \cdot \mathcal{F} := \bigoplus_{k \in \mathbb{Z}} \mathcal{F}^{\oplus a_k} [k].
\end{equation*}
\subsection{Lusztig's geometric categorification of \texorpdfstring{$U^+$}{U+}}\label{section: geometric categorification}
We retain the notation from \cref{section: quantum groups}: $\Phi$ is a simply laced root system with positive roots $\Phi^+$ and simple roots $\{\textbf{e}_i \mid i \in I\}$, $Q$ is a quiver whose underlying graph is the Dynkin diagram of $\Phi$, $R_{\textbf{v}}$ is the space of $\textbf{v}$-dimensional representations of $Q$ and $\GL_{\textbf{v}}$ is the corresponding change of basis group. The $\GL_{\textbf{v}}$-stabilizer of any point in $R_{\textbf{v}}$ is connected. Hence, all equivariant local systems on $\GL _{\textbf{v}}$-orbits in $R_{\textbf{v}}$ are trivial and there is a canonical bijection
\begin{equation*}
\{\GL _{\textbf{v}} \text{-orbits on } R_{\textbf{v}} \}   \overset{1:1} {\leftrightarrow} \Irr(\Perv_G(X)).
\end{equation*}
Combining this with \eqref{eq: kostant partitions parametrize orbits}, we obtain a canonical bijection
\begin{align*}
\KP(\textbf{v}) &\overset{1:1}{\leftrightarrow} \Irr(\Perv(X)) \\
\textbf{c} & \mapsto \IC_{\textbf{c}}.
\end{align*}
This yields a $\mathbb{Z}[v,v^{-1}]$-basis
\begin{equation*}
\{ [\IC_{\textbf{c}}] \mid \textbf{c} \in \KP_{\textbf{v}}\}
\end{equation*}
of $K_{\oplus}(R_{\textbf{v}})$. This is also a $\mathbb{Q}(v)$-basis of
\begin{equation*}
K_{\oplus}(R_{\textbf{v}})_{\mathbb{Q}(v)} :=  \mathbb{Q}(v) \otimes_{\mathbb{Z}[v,v^{-1}]} K_{\oplus}(R_{\textbf{v}}).
\end{equation*}
For any $\textbf{v}, \textbf{w} \in \mathbb{N}[I]$ there is an induction functor \cite[9.2.5]{lusztig2010introduction}
\begin{equation*}
D^b_{\GL_{\textbf{v}} \times \GL_{\textbf{w}}} (R_{\textbf{v}} \times R_{\textbf{w}} ) \overset{\Ind^{\textbf{v} + \textbf{w}}_{\textbf{v}, \textbf{w}}}{\longrightarrow} D^b_{\GL_{\textbf{v} + \textbf{w}}} (R_{\textbf{v} + \textbf{w}})
\end{equation*}
which induces a `convolution map'
\begin{align*}
K_{\oplus}(R_{\textbf{v}}) \otimes_{\mathbb{Z}[v,v^{-1}]} K_{\oplus}(R_{\textbf{w}}) & \rightarrow K_{\oplus}(R_{\textbf{v} + \textbf{w}}) \\
[\mathcal{F}] \otimes [\mathcal{G}] & \mapsto [\Ind^{\textbf{v} + \textbf{w}}_{\textbf{v}, \textbf{w}}( \mathcal{F} \boxtimes \mathcal{G})].
\end{align*}
This equips
\begin{equation*}
K_{\oplus}(Q) := \bigoplus_{\textbf{v} \in \mathbb{N}[I]} K_{\oplus}(R_{\textbf{v}})_{\mathbb{Q}(v)}
\end{equation*}
with the structure of a $\mathbb{Q}(v)$-algebra. For any $\mathcal{F}, \mathcal{G} \in \Semis_{\GL_{\textbf{v}}}(R_{\textbf{v}})$, we define
\begin{equation*}
( \mathcal{F} , \mathcal{G})  = \sum_{k  \in \mathbb{Z}} \dim \Hom^k(\mathbb{D} \mathcal{F},  \mathcal{G}) \cdot v^{-k} \in \mathbb{Q}(v).
\end{equation*}
A priori, this pairing takes values in $\mathbb{Z}((v^{-1}))$ but it can be shown that on semisimple complexes it actually takes values in $\mathbb{Z}((v^{-1})) \cap \mathbb{Q}(v)$. There is an induced bilinear pairing 
\begin{equation*}
(\cdot, \cdot) : K_{\oplus}(R_{\textbf{v}}) \times K_{\oplus}(R_{\textbf{v}}) \rightarrow \mathbb{Q}(v)
\end{equation*}
which can be extended to a symmetric bilinear form
\begin{equation}\label{eq: geometric bilinear form}
K_{\oplus}(Q) \times K_{\oplus}(Q) \rightarrow \mathbb{Q}(v)
\end{equation}
by stipulating that $(K_{\oplus}(R_{\textbf{v}}),K_{\oplus}(R_{\textbf{w}})) = 0$ for any $\textbf{v} \neq \textbf{w}$. The following theorem is Lusztig's geometric categorification of $U^+$.
\begin{theorem}\cite[13.2.11]{lusztig2010introduction}\label{thm: geometric categorification of qua group}
For each $\textbf{v} \in \mathbb{N}[I]$ there is an isomorphism of $\mathbb{Q}(v)$-vector spaces
\begin{equation}\label{eq: Lusztigs categorification iso}
\begin{aligned}
K_{\oplus}(R_{\textbf{v}})_{\mathbb{Q}(v)} & \overset{\sim}{\rightarrow}  U^+_{\textbf{v}} \\
[\IC_{\textbf{c}}] & \mapsto b_{\textbf{c}}.
\end{aligned}
\end{equation} 
This induces an isomorphism of $\mathbb{N}[I]$-graded algebras
\begin{equation*}
K_{\oplus}(Q) \overset{\sim}{\rightarrow}  U^+
\end{equation*}
which intertwines the bilinear forms from \eqref{eq: geometric bilinear form} and \eqref{eq: bilinear form on qua groups}.
\end{theorem}
The change of basis matrix between $\Bstdv$ and $\Bcanv$ also has a geometric interpretation.
\begin{proposition}\label{prop: geometric interpretation of base change between standard and canonical basis}\cite[10.7]{lusztig1990canonical}
The entries of the change of basis matrix $P$ between $\Bstdv$ and $\Bcanv$ are given by
\begin{equation*}
p_{\textbf{c}, \textbf{c}'}(v) = \sum_{k \in \mathbb{Z}} \dim H^k( \iota_x^* \IC_{\textbf{c}}) v^k
\end{equation*}
where $\iota_x : \{x \} \hookrightarrow \mathcal{O}_{\textbf{c}'}$ is any base point.
\end{proposition}
\begin{remark}
When $Q$ is the equioriented type A quiver, the polynomials $p_{\textbf{c}, \textbf{c}'}(v)$ turn out to be certain parabolic Kazhdan-Lusztig polynomials \cite{zelevinskii1985two}.
\end{remark}
Next, we discuss the geometric interpretation of monomials. Let
\begin{align*}
\textbf{i} &= (i_1, ..., i_m) \in I^m \\
\textbf{a} &= (a_1, ..., a_m) \in \mathbb{N}^m
\end{align*}
such that $(\textbf{i}, \textbf{a})$ has weight $\textbf{v}$. Let $B \subset \GL_{\textbf{v}}$ be a fixed Borel subgroup. Then there is a unique $B$-stable flag $F^*$ of $I$-graded subspaces 
\begin{equation*}
\mathbb{C}^{\textbf{v}} = F^0 \supset F^1 \supset ... \supset F^m = \{0\}
\end{equation*}
of type $(\textbf{i}, \textbf{a})$ (i.e. $F^{k-1}/F^k$ is pure of weight $i_k$ and dimension $a_k$). The stabilizer $P_{\textbf{i}, \textbf{a}} \subset \GL_{\textbf{v}}$ of $F^*$ is a parabolic subgroup containing $B$. Let $R_{\textbf{i}, \textbf{a}} \subset R_{\textbf{v}}$ be the subspace of all representations $x$ that are compatible with the flag $F^*$, i.e. $x(F^k) \subset F^k$ for all $k$. We define
\begin{equation*}
\tilde{R}_{\textbf{i}, \textbf{a}} := \GL_{\textbf{v}} \times^{P_{\textbf{i}, \textbf{a}}} R_{\textbf{i}, \textbf{a}}.
\end{equation*}
Then the morphism
\begin{align*}
\mu_{\textbf{i}, \textbf{a}} : \tilde{R}_{\textbf{i}, \textbf{a}} & \rightarrow R_{\textbf{v}}\\
(g,v) & \mapsto g\cdot v
\end{align*}
is proper and we define
\begin{equation*}
\Sbf_{\textbf{i}, \textbf{a}} := (\mu_{\textbf{i}, \textbf{a}})_* \mathcal{C}_{\tilde{R}_{\textbf{i}, \textbf{a}}}.
\end{equation*}
By the decomposition theorem \cite{beilinson2018faisceaux}, we have
\begin{equation*}
\Sbf_{\textbf{i}, \textbf{a}}  \in \Semis_{\GL_{\textbf{v}}}(R_{\textbf{v}}).
\end{equation*}
Moreover,
\begin{equation}\label{eq: springer sheaves are verdier invariant}
\begin{aligned}
\mathbb{D} \Sbf_{\textbf{i}, \textbf{a}} & = \mathbb{D} (\mu_{\textbf{i}, \textbf{a}})_* \mathcal{C}_{\tilde{R}_{\textbf{i}, \textbf{a}}} \\
&  =  (\mu_{\textbf{i}, \textbf{a}})_* \mathbb{D} \mathcal{C}_{\tilde{R}_{\textbf{i}, \textbf{a}}} \\
& =  (\mu_{\textbf{i}, \textbf{a}})_*  \mathcal{C}_{\tilde{R}_{\textbf{i}, \textbf{a}}}  \\
& =  \Sbf_{\textbf{i}, \textbf{a}}.
\end{aligned}
\end{equation}
\begin{proposition}\label{prop: i a springer sheaf identifies with monomial}
The isomorphism $K_{\oplus}(R_{\textbf{v}})_{\mathbb{Q}(v)} \overset{\sim}{\rightarrow}  U^+_{\textbf{v}} $ from \eqref{eq: Lusztigs categorification iso} identifies $[\Sbf_{\textbf{i}, \textbf{a}}]$ with $E_{\textbf{i}, \textbf{a}}$.
\end{proposition}
\begin{proof}
See \cite[\S 5]{reineke2003quivers} for an argument in the language of Hall algebras. In the language of perverse sheaves this can be explained as follows: It is straightforward to check that $E_{i_j}^{(a_j)}$ is a canonical basis element for all $j=1,...,m$. Hence, we have
\begin{equation*}
E_{\textbf{i}, \textbf{a}} =   E_{i_1}^{(a_1)} E_{i_2}^{(a_2)} ... E_{i_m}^{(a_m)} \overset{\eqref{eq: Lusztigs categorification iso}}{=} [ \IC_{\textbf{e}_{i_1}^{a_1}}] \cdot [ \IC_{\textbf{e}_{i_2}^{a_2}}] \cdot ... \cdot [ \IC_{\textbf{e}_{i_m}^{a_m}}]
\end{equation*}
Note that $R_{\textbf{e}_{i_j}^{a_j}} = \{ 0 \}$ is just a point and hence $\IC_{\textbf{e}_{i_j}^{a_j}}$ is just the constant sheaf $\mathcal{C}_{R_{\textbf{e}_{i_j}^{a_j}}}$. Spelling out the definition of the corresponding convolution explicitly (see for example \cite[(10.3.8)]{achar2021perverse}), we arrive at
\begin{equation*}
E_{\textbf{i}, \textbf{a}} = [ \mathcal{C}_{R_{\textbf{e}_{i_1}^{a_1}}}] \cdot [ \mathcal{C}_{R_{\textbf{e}_{i_2}^{a_2}}}] \cdot ... \cdot [ \mathcal{C}_{R_{\textbf{e}_{i_m}^{a_m}}}] =  [(\mu_{\textbf{i}, \textbf{a}})_* \mathcal{C}_{\tilde{R}_{\textbf{i}, \textbf{a}}}] = [ \Sbf_{\textbf{i}, \textbf{a}}].
\end{equation*}
\end{proof}
Combining \cref{thm: geometric categorification of qua group} and \cref{prop: i a springer sheaf identifies with monomial} yields the following geometric interpretation of the matrix $\Psi$ from \eqref{eq: def of psi in qua group case}:
\begin{equation}\label{eq: geometric description of inner product of monomials}
\begin{aligned}
\Psi_{\textbf{c}, \textbf{c}'} & = (E_{\textbf{i}, \textbf{a}_{\textbf{c}}},  E_{\textbf{i}, \textbf{a}_{\textbf{c}'}}) \\
& = ([\Sbf_{\textbf{i}, \textbf{a}_{\textbf{c}}}],  [\Sbf_{\textbf{i}, \textbf{a}_{\textbf{c}'}}]) \\
& = \sum_{k \in \mathbb{Z}} \dim \Hom^k (\mathbb{D} (\Sbf_{\textbf{i}, \textbf{a}_{\textbf{c}}}),  \Sbf_{\textbf{i}, \textbf{a}_{\textbf{c}'}}) \cdot v^{-k} \\
\overset{\eqref{eq: springer sheaves are verdier invariant}}&{=} \sum_{k \in \mathbb{Z}} \dim \Hom^k(\Sbf_{\textbf{i}, \textbf{a}_{\textbf{c}}} ,  \Sbf_{\textbf{i}, \textbf{a}_{\textbf{c}'}}) \cdot v^{-k}.
\end{aligned}
\end{equation}
There also is a geometric version of \cref{lem: existence of monomial basis}.
\begin{lemma}\cite[2.2]{reineke2003quivers}\label{lem: existence of res of singularities}
We can find $m \in \mathbb{N}$, $\textbf{i} \in I^m$, and for each $\textbf{c} \in \KP(\textbf{v})$ a sequence $\textbf{a}_{\textbf{c}} \in \mathbb{N}^m$ such that $(\textbf{i}, \textbf{a}_{\textbf{c}})$ has weight $\textbf{v}$ and
\begin{enumerate}
\item $\mu_{\textbf{i}, \textbf{a}_{\textbf{c}}}$ has image $\overline{\mathcal{O}}_{\textbf{c}}$,
\item $\mu_{\textbf{i}, \textbf{a}_{\textbf{c}}}$ restricts to an isomorphism $\mu_{\textbf{i}, \textbf{a}_{\textbf{c}}}^{-1}(\mathcal{O}_{\textbf{c}}) \overset{\sim}{\rightarrow} \mathcal{O}_{\textbf{c}}$ over $\mathcal{O}_{\textbf{c}}$,
\end{enumerate}
for all $\textbf{c} \in \KP(\textbf{v})$.
\end{lemma}
Note that the lemma above (together with the decomposition theorem) implies that 
\begin{equation*}
\Sbf_{\textbf{i}, \textbf{a}_{\textbf{c}}} \cong \IC_{\textbf{c}} \oplus \bigoplus_{k=1}^n \IC_{\textbf{c}_k} [d_k]
\end{equation*}
for some (not necessarily distinct) Kostant partitions $\textbf{c}_k < \textbf{c}$ and  integers $d_k \in \mathbb{Z}$. Taking the corresponding classes in $K_{\oplus}(R_{\textbf{v}})_{\mathbb{Q}(v)} \cong  U^+_{\textbf{v}}$  yields
\begin{equation*}
E_{\textbf{i}, \textbf{a}_{\textbf{c}}} = [\Sbf_{\textbf{i}, \textbf{a}_{\textbf{c}}}] = [\IC_{\textbf{c}}] +  \sum_{k=1}^n [\IC_{\textbf{c}_k} [d_k] ] = E_{\textbf{c}} + \sum_{k=1}^n v^{d_k} E_{\textbf{c}_k}.
\end{equation*}
Hence, \cref{lem: existence of monomial basis} and the existence of a monomial basis is a consequence of \cref{lem: existence of res of singularities}.
\subsection{Computation of some Ext groups}
In this section we derive a formula for the graded dimensions of $\Hom$-spaces between certain semisimple complexes. Let $G$ be a connected complex reductive group with a fixed Borel subgroup $B \subset G$, torus $T\subset B$ and Weyl group $W$. Let $P_1, P_2 \subset G$ be two parabolic subgroups containing $B$ with corresponding parabolic subgroups $W_1,W_2 \subset W$. For any $w \in W$, we write ${}^wP_i := \dot{w} P_i \dot{w}^{-1}$ ($i=1,2$) where $\dot{w} \in N_G(T)$ is a lift of $w$. Consider the orbit partition
\begin{equation}\label{Orbit partition on partial flag variety product}
G/P_1 \times G/P_2 = \bigsqcup_{[w] \in W_1\backslash W /W_2} Y_{[w]}
\end{equation}
where $Y_{[w]} = G \cdot (P_1,{}^wP_2) \cong G/(P_1 \cap {}^wP_2)$. We start by computing the Poincar\'e polynomial of $Y_{[w]}$. 
\begin{lemma}\label{Computation of Yw}
We have
\begin{equation*}
\chi^G_v(Y_{[w]}) = v^{-r_{[w]}} \frac{\chi^G_v(G/B)}{\chi_v(P_1/B)\chi_v(P_2/B)} \sum_{y \in [w]}v^{2l(y)}
\end{equation*}
where $r_{[w]} = \dim Y_{[w]} + \dim P_1/B + \dim P_2/B - \dim G/B$
\end{lemma}
\begin{proof}
Let
\begin{equation*}
q: G/B \times G/B \rightarrow G/P_1 \times G/P_2
\end{equation*}
be the projection and
\begin{equation*}
X_{[w]} := q^{-1}(Y_{[w]}).
\end{equation*}
We first compute $\chi^G_v(X_{[w]})$. Note that $X_{[w]} = \bigsqcup_{y \in [w]} \mathcal{O}_y$ where
\begin{equation*}
G/B \times G/B = \bigsqcup_{y \in W} \mathcal{O}_y
\end{equation*}
is the partition into $G$-orbits. For any $y\in W$, the map $pr_1: \mathcal{O}_y \rightarrow G/B$ is a $G$-equivariant locally trivial fibration with fiber $\mathbb{A}^{\dim \mathcal{O}_y - \dim G/B}$. Hence, $H^k_G(\mathcal{O}_y) \cong H^k_G(G/B)$ for all $k \ge 0$ and
\begin{equation}\label{eq: poincare poly of Oy and G/B}
\chi^G_v(\mathcal{O}_y) = v^{ \dim \mathcal{O}_y  - \dim G/B }\chi^G_v(G/B).
\end{equation}
Moreover,
\begin{equation*}
H^G_k(\mathcal{O}_y)\overset{PD}{\cong} H_G^{2 \dim \mathcal{O}_y - k}(\mathcal{O}_y) \cong  H_G^{2 \dim \mathcal{O}_y - k}(G/B)
\end{equation*}
vanishes for $k$ odd. By \cref{parity vanishing direct sum of Borel Moore}, this implies
\begin{equation}\label{BM orbit decomposition}
H^G_k(X_{[w]}) = \bigoplus_{y \in [w]} H^G_k(\mathcal{O}_y).
\end{equation}
Hence, we get
\begin{equation}\label{Computation of Xw}
\begin{aligned}
\chi^G_v(X_{[w]}) &\overset{\eqref{eq: def of poincare poly}}{=} v^{\dim X_{[w]}} \sum_{k\in \mathbb{Z}} \dim H_G^k( X_{[w]})v^{-k} \\
& \overset{PD}{=} v^{\dim X_{[w]}} \sum_{k\in \mathbb{Z}} \dim H^G_{2\dim X_{[w]} - k}( X_{[w]})v^{-k} \\
& = v^{-\dim X_{[w]}} \sum_{k\in \mathbb{Z}}   \dim H^G_{k}( X_{[w]})v^{k} \\
& \overset{\eqref{BM orbit decomposition}}{=} v^{-\dim X_{[w]}} \sum_{k\in \mathbb{Z}} \sum_{y \in [w]} \dim H^G_k(\mathcal{O}_y)v^{k} \\
&=  v^{-\dim X_{[w]}} \sum_{k\in \mathbb{Z}} \sum_{y \in [w]}   \dim H^G_{2\dim\mathcal{O}_y- k}(\mathcal{O}_y) v^{2\dim\mathcal{O}_y-k}\\
& \overset{PD}{=} v^{-\dim X_{[w]}} \sum_{k\in \mathbb{Z}} \sum_{y \in [w]} \dim H_G^k(\mathcal{O}_y) v^{2\dim\mathcal{O}_y-k} \\
&= v^{-\dim X_{[w]}} \sum_{y \in [w]} \chi^G_v(\mathcal{O}_y) v^{\dim\mathcal{O}_y} \\
&\overset{\eqref{eq: poincare poly of Oy and G/B}}{=} v^{-\dim X_{[w]}} \sum_{y \in [w]} \chi^G_v(G/B) v^{2 \dim\mathcal{O}_y-\dim G/B } .
\end{aligned}
\end{equation}
Note that $ \dim\mathcal{O}_y = \dim G/B + l(y)$, so we can also write \eqref{Computation of Xw} as 
\begin{equation}\label{alternative computation of Xw}
\begin{aligned}
\chi^G_v(X_{[w]}) & = v^{-\dim X_{[w]}} \sum_{y \in [w]} \chi^G_v(G/B) v^{2 (\dim G/B + l(y) ) - \dim G/B } \\
&= v^{\dim G/B-\dim X_{[w]}} \chi^G_v(G/B)\sum_{y \in [w]} v^{2l(y)}.
\end{aligned}
\end{equation}
Let $\tilde{q} : X_{[w]} \rightarrow Y_{[w]}$ be the restriction of the map $q$ to $Y_{[w]}$. This is a $G$-equivariant locally trivial fibration with fiber $P_1/B \times P_2/B$. Hence,
\begin{equation}\label{computation of t exponent}
\dim X_{[w]} -\dim G/B = \dim Y_{[w]} + \dim P_1/B + \dim P_2/B - \dim G/B = r_{[w]}.
\end{equation}
To complete the proof, we claim that
\begin{equation}\label{chi(X) and chi(Y)}
\chi^G_v(X_{[w]}) = \chi^G_v(Y_{[w]})\chi_v(P_1/B)\chi_v(P_2/B).
\end{equation}
Note that $\tilde{q}$ is a $G$-equivariant proper morphisms, so by the decomposition theorem $\tilde{q}_* \mathcal{C}_{X_{[w]}}$ is a semisimple $G$-equivariant complex on $Y_{[w]}$. The only simple $G$-equivariant perverse sheaf on $Y_{[w]}$ is $\mathcal{C}_{Y_{[w]}}$. Therefore, we can write
\begin{equation}\label{pushforward of constant sheaf on Xw}
\tilde{q}_* \mathcal{C}_{X_{[w]}} = f(v) \cdot \mathcal{C}_{Y_{[w]}}
\end{equation}
for some $f(v) \in \mathbb{N}[v, v^{-1}]$. Pushing \eqref{pushforward of constant sheaf on Xw} forward to the point, we get from \eqref{eq: Groth group interpretation of poincare poly} that
\begin{equation}\label{eq: chi of Yw up to f(v)}
\chi^G_v(X_{[w]}) =  f(v) \cdot \chi^G_v(Y_{[w]}).
\end{equation}
Consider the inclusion maps $\iota_x : \{x\} \hookrightarrow \mathcal{O}_y$ and $\iota_{q^{-1}(x)}: q^{-1}(x) \hookrightarrow X_{[w]}$. Using proper base change (after forgetting the $G$-equivariance), we get
\begin{align*}
f(v) \cdot \mathcal{C}_{\{x\}}  & = f(v)  \cdot \iota_x^* \mathcal{C}_{Y_{[w]}} [- \dim Y_{[w]}] \\
\overset{\eqref{pushforward of constant sheaf on Xw}}&{=}  \iota_x^*\tilde{q}_* \mathcal{C}_{X_{[w]}} [- \dim Y_{[w]}] \\
& = \tilde{q}_* \iota_{q^{-1}(x)}^* \mathcal{C}_{X_{[w]}} [- \dim Y_{[w]}] \\
& = \tilde{q}_* \mathcal{C}_{q^{-1}(x)}.
\end{align*}
In the Grothendieck group this implies
\begin{equation*}
f(v) \cdot [ \mathcal{C}_{\{x\}}] = [ \tilde{q}_* \mathcal{C}_{q^{-1}(x)}] \overset{\eqref{eq: Groth group interpretation of poincare poly}}{=} \chi_v(q^{-1}(x))[ \mathcal{C}_{\{x\}}].
\end{equation*}
Since $q^{-1}(x) \cong P_1/B \times P_2/B$, we get
\begin{equation*}
f(v) =  \chi_v(q^{-1}(x)) = \chi_v(P_1/B)\chi_v(P_2/B).
\end{equation*}
Thus \eqref{chi(X) and chi(Y)} follows from \eqref{eq: chi of Yw up to f(v)}. Hence, we obtain
\begin{align*}
\chi^G_v(Y_{[w]}) \overset{\eqref{chi(X) and chi(Y)}}&{=} \frac{1}{\chi_v(P_1/B)\chi_v(P_2/B)} \chi^G_v(X_{[w]}) \\
 \overset{\eqref{alternative computation of Xw},\eqref{computation of t exponent}}&{=} v^{-r_{[w]}} \frac{\chi^G_v(G/B)}{\chi_v(P_1/B)\chi_v(P_2/B)} \sum_{y \in [w]}v^{2l(y)} .
\end{align*}
\end{proof}
We now use the lemma above to compute the graded dimensions of certain $\Hom$-spaces which includes the entries of the matrix $\Psi$ from \eqref{eq: def of psi in qua group case} as a special case. Let $V$ be a $G$-representation and $V_i\subset V$ a $P_i$-stable subspace for $i=1,2$. We can then consider the $G$-space
\begin{equation}\label{resoltion of Springer type}
\tilde{V}_i:= G \times^{P_i} V_i = \{(gP_i,v) \in G/P_i \times V \mid v \in gV_i \}.
\end{equation}
This comes with the two projection maps
\begin{equation*}
\begin{tikzcd}
            & \tilde{V}_i \arrow[ld, "\pi_i"'] \arrow[rd, "\mu_i"] &   \\
G/P_i &                                                      & V
\end{tikzcd}
\end{equation*}
and the associated `Springer sheaf'
\begin{equation*}
\Sbf_i:= (\mu_i)_* \mathcal{C}_{\tilde{V}_i} \in D^b_G(V).
\end{equation*}
For any $w \in W$ we can define the subspace
\begin{equation*}
{}^w V_i := \{ \dot{w}\cdot v \mid v \in V_i\} \subset V
\end{equation*}
where $\dot{w} \in N_G(T)$ is a lift of $w \in W$.
\begin{theorem}\label{theorem: Ext formula}
We have
\begin{equation*}
\sum_{k \in \mathbb{Z}} \dim \Hom^k(\Sbf_1, \Sbf_2) v^{-k} = v^{-\dim  G/B - \dim V_1 - \dim V_2} \frac{\chi^G_v(G/B)}{\chi_v (P_1/B)\chi_v (P_2/B)} \sum_{w \in W} v^{2(l(w)+ \dim V_1\cap {}^w V_2)}.
\end{equation*}
\end{theorem}
\begin{proof}
Consider the Steinberg variety
\begin{equation*}
Z := \tilde{V}_1 \times_V \tilde{V}_2.
\end{equation*}
By \cite[8.6]{chriss2009representation}, there is an isomorphism
\begin{equation}\label{Ext isomorphic to BM homomology}
\Hom^k(\Sbf_1, \Sbf_2) \cong H^G_{d_1+d_2-k}(Z)
\end{equation}
where
\begin{equation*}
d_i = \dim \tilde{V}_i = \dim G - \dim P_i + \dim V_i.
\end{equation*}
The orbit partition on $G/P_1 \times G/P_2$ from \eqref{Orbit partition on partial flag variety product} induces a partition
\begin{equation*}
Z = \bigsqcup_{[w]\in W_1 \backslash W / W_2} Z_{[w]}
\end{equation*}
where $Z_{[w]} := (\pi_1 \times \pi_2)^{-1} (Y_{[w]})$. The projection $Z_{[w]} \rightarrow Y_{[w]}$ is a $G$-equivariant vector bundle with fiber $V_1 \cap {}^w V_2$. Hence, $Z_{[w]}$ is smooth and
\begin{equation}\label{PD for parts of Z}
H^G_k(Z_{[w]}) \overset{PD}{=}H_G^{2d_{[w]}-k}(Z_{[w]}) \cong H_G^{2d_{[w]}-k}(Y_{[w]}).
\end{equation}
where
\begin{equation*}
d_{[w]} := \dim Z_{[w]} = \dim Y_{[w]} + \dim V_1 \cap {}^w V_2.
\end{equation*}
It follows from \cref{Computation of Yw} that $H_G^{k}(Y_{[w]})$ vanishes for $k$ odd, so $H^G_k(Z_{[w]})$ also vanishes for $k$ odd. By \cref{parity vanishing direct sum of Borel Moore}, we get
\begin{equation}\label{decomposition of Z by orbits}
H^G_k(Z) = \bigoplus_{[w] \in W_1 \backslash W /W_2} H^G_k(Z_{[w]}).
\end{equation}
Hence, we have
\begin{align*}
(\dag) & := \sum_{k \in \mathbb{Z}}  \dim \Hom^k(\Sbf_1, \Sbf_2) v^{-k}  \\
\overset{\eqref{Ext isomorphic to BM homomology}}&{\cong} \sum_{k \in \mathbb{Z}}  \dim H^G_{d_1+d_2-k}(Z)v^{-k}\\
\overset{\eqref{decomposition of Z by orbits}}&{=} \sum_{k \in \mathbb{Z}} \sum_{[w] \in W_1 \backslash W /W_2}  \dim H^G_{d_1 + d_2 -k} (Z_{[w]})v^{-k} \\
\overset{\eqref{PD for parts of Z}}&{=}  \sum_{k \in \mathbb{Z}} \sum_{[w] \in W_1 \backslash W /W_2}  \dim H_G^{2d_{[w]}-d_1 -d_2 + k}(Y_{[w]})v^{-k} \\
&=  \sum_{k \in \mathbb{Z}} \sum_{[w] \in W_1 \backslash W /W_2} v^{2d_{[w]}-d_1-d_2} \dim H_G^{k}(Y_{[w]})v^{-k} \\
\overset{\eqref{eq: def of poincare poly}}&{=} \sum_{[w] \in W_1 \backslash W /W_2} v^{2d_{[w]}- d_1-d_2 - \dim Y_{[w]}}\chi^G_{v}(Y_{[w]})\\
\overset{\ref{Computation of Yw}}&{=}  \sum_{[w] \in W_1 \backslash W /W_2} v^{2d_{[w]} -d_1-d_2- \dim Y_{[w]}- r_{[w]} } \frac{\chi^G_v(G/B)}{\chi_v(P_1/B)\chi_v(P_2/B)} \sum_{y \in [w]}v^{2l(y)} \\
&=\frac{\chi^G_v(G/B)}{\chi_v(P_1/B)\chi_v(P_2/B)}  \sum_{w \in W} v^{2d_{[w]} -d_1 -d_2- \dim Y_{[w]}-r_{[w]} +2l(w)}.
\end{align*}
Recall that
\begin{align*}
d_{[w]} &= \dim Y_{[w]} + \dim V_1\cap {}^w V_2 \\
r_{[w]} &= \dim Y_{[w]} + \dim P_1/B + \dim P_2/B - \dim G/B \\
d_i &= \dim G - \dim P_i + \dim V_i.
\end{align*}
Hence, we have
\begin{equation*}
2d_{[w]} -d_1 -d_2- \dim Y_{[w]}-r_{[w]}  =  2\dim V_1\cap {}^w V_2 - \dim G/B-\dim V_1 - \dim V_2
\end{equation*}
and thus
\begin{align*}
(\dag)&=\frac{\chi_v(G/B)}{\chi_v(P_1/B) \chi_v(P_2/B)}  \sum_{w \in W} v^{2\dim V_1 \cap {}^w V_2 - \dim G/B-\dim V_1 - \dim V_2 + 2l(w)} \\
&= v^{-\dim G/B- \dim V_1 - \dim V_2} \frac{\chi_v(G/B)}{\chi_v (P_1/B)\chi_v (P_2/B)} \sum_{w \in W} v^{2(l(w)+ \dim V_1\cap wV_2)}.
\end{align*}
\end{proof}
\begin{remark}
The theorem above can be used to compute the graded dimension of some interesting Ext-algebras. For example, KLR algebras can be realised as Ext algebras of a direct sum of sheaves of the form $\Sbf_i$ as above \cite{varagnolo2011canonical}. Similarly, affine Hecke algebras at a central character (away from roots of unity) also have such a realization \cite{chriss2009representation}. We will make use of this realization in \cref{Section: dimensions of simples} to compute the dimensions of the simple modules for the affine Hecke algebra of $\GL_n$. In both of these examples, the parabolic subgroups $P_i$ in the definition of the corresponding sheaves $\Sbf_i$ are always Borel subgroups. However, to compute the bilinear form on $U^+$ at monomials, we also need to consider parabolic subgroups that are not Borel subgroups.
\end{remark}
We now apply the formula from the theorem above to the quantum groups setting. Recall for this that $S_{\textbf{v}} = \prod_{i \in I} S_{\textbf{v}_i}$ is the Weyl group of $\GL_{\textbf{v}} = \prod_{i \in I} \GL_{\textbf{v}_i}$, $B \subset \GL_{\textbf{v}}$ is fixed a Borel subgroup and $R_{\textbf{i}, \textbf{a}} \subset R_{\textbf{v}}$ is the subspace that fixes the unique $B$-stable flag $F^*$ of type $(\textbf{i}, \textbf{a})$. For any sequence $\textbf{a} = (a_1,a_2, ..., a_m) \in \mathbb{N}^m$ we define
\begin{equation*}
[\textbf{a}]! := [a_1]! [a_2]! \cdot ... \cdot [a_m]!.
\end{equation*}
\begin{corollary}\label{cor: geometric counting formula for inner product}
Let $(\textbf{i}, \textbf{a})$ and $(\textbf{i}', \textbf{a}')$ be of weight $\textbf{v}$. Then
\begin{equation*}
(E_{\textbf{i}, \textbf{a}}, E_{\textbf{i}', \textbf{a}'}) = v^{- \dim R_{\textbf{i}, \textbf{a}} - \dim R_{ \textbf{i}', \textbf{a}'}}  \frac{\prod_{i \in I} (1-v^{-2})^{-\textbf{v}_i}}{ [\textbf{a}]! [\textbf{a}']!} \sum_{w \in S_{\textbf{v}}} v^{2(l(w) + \dim R_{\textbf{i}, \textbf{a}} \cap  {}^w R_{\textbf{i}', \textbf{a}'})}.
\end{equation*}
\end{corollary}
\begin{proof}
The parabolics $P_{\textbf{i}, \textbf{a}}, P_{\textbf{i}, \textbf{a}'}$ contain the common Borel subgroup $B$. Hence, we can apply \cref{theorem: Ext formula} to obtain
\begin{align*}
(E_{\textbf{i}, \textbf{a}}, E_{\textbf{i}', \textbf{a}'}) & = \sum_{k \in \mathbb{Z}} \Hom^k(\Sbf_{\textbf{i}, \textbf{a}}, \Sbf_{\textbf{i}', \textbf{a}'} ) v^{-k} \\
&= v^{-s} \frac{\chi^{\GL_{\textbf{v}}}_v(\GL_{\textbf{v}}/B)}{\chi_v (P_{\textbf{i}, \textbf{a}}/B)\chi_v (P_{\textbf{i}', \textbf{a}'}/B)} \sum_{w \in S_{\textbf{v}}} v^{2(l(w) + \dim R_{\textbf{i}, \textbf{a}} \cap {}^w R_{\textbf{i}', \textbf{a}'})}
\end{align*}
where
\begin{equation*}
s = \dim \GL_{\textbf{v}}/B + \dim R_{\textbf{i}, \textbf{a}} + \dim R_{\textbf{i}, \textbf{a}}.
\end{equation*}
Note that the Levi factor of $P_{\textbf{i}, \textbf{a}}$ is isomorphic to $\GL_{\textbf{a}_1} \times ... \times \GL_{\textbf{a}_m}$. The Poincar\'e polynomials for the flag variety are well-known:
\begin{align*}
\chi^{\GL_{\textbf{v}}}_v(\GL_{\textbf{v}}/B) & = v^{\dim \GL_{\textbf{v}}/B} \chi^B_v(\pt) = v^{\dim \GL_{\textbf{v}}/B}\prod_{i \in I} (1-v^{-2})^{-\textbf{v}_i}  \\
\chi_v (P_{\textbf{i}, \textbf{a}}/B) &= [\textbf{a}]! \\
\chi_v (P_{\textbf{i}', \textbf{a}'}/B) &= [\textbf{a}']!.
\end{align*} 
Thus we have
\begin{align*}
(E_{\textbf{i}, \textbf{a}}, E_{\textbf{i}', \textbf{a}'}) &= v^{- s}     \frac{ v^{\dim \GL_{\textbf{v}}/B} \prod_{i \in I}   (1-v^{-2})^{-\textbf{v}_i} }{[\textbf{a}]! [\textbf{a}']!} \sum_{w \in S_{\textbf{v}}} v^{2(l(w) + \dim R_{\textbf{i}, \textbf{a}} \cap  {}^w R_{\textbf{i}', \textbf{a}'})} \\
&=     v^{- \dim R_{\textbf{i}, \textbf{a}} - \dim R_{ \textbf{i}', \textbf{a}'}}  \frac{\prod_{i \in I} (1-v^{-2})^{-\textbf{v}_i}}{ [\textbf{a}]! [\textbf{a}']!} \sum_{w \in S_{\textbf{v}}} v^{2(l(w) + \dim R_{\textbf{i}, \textbf{a}} \cap  {}^w R_{\textbf{i}', \textbf{a}'})}.
\end{align*}
\end{proof}
\section{Examples in type A}\label{section: examples}
In order to be able to execute \cref{algorithm: new algo}, we need to explicitly compute
\begin{itemize}
\item The sequences $\textbf{i}$ and $\textbf{a}_c$ from \cref{lem: existence of res of singularities} that give rise to the monomial basis
\begin{equation*}
\Bmon = \{ E_{\textbf{i}, \textbf{a}_{\textbf{c}}}  \mid \textbf{c} \in \KP(\textbf{v}) \}.
\end{equation*}
\item The partial order $\le$ on $\KP(\textbf{v})$, i.e. the dimensions of the orbits $\mathcal{O}_{\textbf{c}}$;
\end{itemize}
The orbit dimensions can be computed using Lusztig's formula \cite[Proposition 6.6]{lusztig1990canonical}. Alternatively, they can also be computed from the sequences $\textbf{i}$ and $\textbf{a}_c$ using that
\begin{equation}\label{eq: orbit dimension formula}
\dim \mathcal{O}_{\textbf{c}} = \dim \tilde{R}_{\textbf{i}, \textbf{a}_{\textbf{c}}} = \dim G - \dim P_{\textbf{i}, \textbf{a}_{\textbf{c}}} + R_{\textbf{i}, \textbf{a}_{\textbf{c}}}
\end{equation}
by \cref{lem: existence of res of singularities}. We now recall an explicit description of the sequences $\textbf{i}, \textbf{a}_{\textbf{c}}$ from \cite{reineke2003quivers} for the equioriented type A quiver
\begin{equation}\label{eq: type A quiver}
A_n = \begin{tikzcd}
\bullet \arrow[r] & \bullet \arrow[r] & ... \arrow[r] & \bullet
\end{tikzcd}
\end{equation}
with $n$ vertices and $n-1$ arrows (see also \cite[\S 4]{achar2018combinatorics} for an alternative combinatorial description using triangular arrays). Let $U^+ = U^+(A_n)$ be the associated positive quantum group. The set of positive roots is given by
\begin{equation*}
\Phi^+ = \{ \alpha_{ij} = \textbf{e}_i + \textbf{e}_{i+1} + ... + \textbf{e}_j \mid i \le j\}
\end{equation*}
where $\textbf{e}_1 , ..., \textbf{e}_n$ are the simple roots labelling the quiver from left to right. The reduced expression
\begin{equation*}
w_0 = s_n (s_{n-1} s_n) ... (s_i s_{i+1} ... s_n) ... (s_1 s_2 ... s_n)
\end{equation*}
is adapted to $A_n$ in the sense of \cite[4.7]{lusztig1990canonical}. Let
\begin{equation*}
\textbf{i} = (n, n-1, n, ...,i, i+1, ..., n, ... 1, 2, ..., n).
\end{equation*}
For any Kostant partition $\textbf{c} = (m_{ij})_{i \le j} \in \KP(\textbf{v})$ set
\begin{equation*}
\textbf{a}_{\textbf{c}} = (m_{nn}, m_{n-1,n-1} + m_{n-1, n} , m_{n-1,n}, ..., m_{1,1} + m_{1,2} + ... + m_{1,n} , ..., m_{1,n}).
\end{equation*}
Then by \cite{reineke2003quivers} the $\textbf{a}_{\textbf{c}}$ satisfy the conditions from \cref{lem: existence of res of singularities}. Hence,
\begin{equation*}
\Bmonv = \{ E_{\textbf{i}, \textbf{a}_{\textbf{c}}} \mid \textbf{c} \in \KP(\textbf{v}) \}
\end{equation*}
is a monomial basis of $U^+_{\textbf{v}}$. We now compute the canonical basis of $U^+_{\textbf{v}}$ in two examples.
\subsection{Example: \texorpdfstring{$\textbf{v} = (2,2)$}{v=(2,2)}}
Let $n=2$ and $\textbf{v} = (2,2)$. The order on the set of positive roots is given by
\begin{equation*}
\Phi^+  = \{ \textbf{e}_2 <  \textbf{e}_1 +  \textbf{e}_2 <  \textbf{e}_1 \}.
\end{equation*}
There are three Kostant partitions of $\textbf{v}$:
\begin{equation*}
\textbf{c}_1 =  ( 2,0,2) ,\quad \textbf{c}_2 = (1,1,1) ,\quad \textbf{c}_3 =   (0,2,0).
\end{equation*}
We have
\begin{equation*}
\textbf{i} = (2,1,2)
\end{equation*}
and
\begin{equation*}
\textbf{a}_{\textbf{c}_1} = (2,2,0), \quad \textbf{a}_{\textbf{c}_2} = (1,2,1), \quad  \textbf{a}_{\textbf{c}_3} = (0,2,2).
\end{equation*}
Let $B \subset \GL_2$ be the standard (geometric) Borel subgroup consisting of lower triangular matrices and consider the corresponding Borel $B \times B \subset \GL_2 \times \GL_2 = \GL_{\textbf{v}}$. The space of quiver representations $R_{\textbf{v}}$ can be identified with the set of $2 \times 2$ matrices $ \begin{spm} * & * \\ * & * \end{spm}$. The corresponding spaces $R_{\textbf{i}, \textbf{a}_{\textbf{c}}} \subset R_{\textbf{v}}$ are given by 
\begin{equation*}
R_{\textbf{i}, \textbf{a}_{\textbf{c}_1}} = \begin{spm} 0 & 0 \\ 0 & 0 \end{spm}, \quad 
R_{\textbf{i}, \textbf{a}_{\textbf{c}_2}} = \begin{spm} 0 & 0 \\ * & * \end{spm}, \quad 
R_{\textbf{i},  \textbf{a}_{\textbf{c}_3}} = \begin{spm} * & * \\ * & * \end{spm} .
\end{equation*}
The corresponding parabolic subgroups of $\GL_{\textbf{v}}$ are
\begin{equation*}
P_{\textbf{i}, \textbf{a}_{\textbf{c}_1}} = \GL_{\textbf{v}}, \quad 
P_{\textbf{i}, \textbf{a}_{\textbf{c}_2}} = \GL_2 \times B, \quad 
P_{\textbf{i},  \textbf{a}_{\textbf{c}_3}} = \GL_{\textbf{v}}.
\end{equation*}
Using \eqref{eq: orbit dimension formula}, we can compute the orbit dimensions
\begin{equation*}
\dim \mathcal{O}_{\textbf{c}_1} = 0, \quad 
\dim \mathcal{O}_{\textbf{c}_2} = 3, \quad 
\dim \mathcal{O}_{\textbf{c}_3} = 4
\end{equation*}
and hence the order on Kostant partitions is
\begin{equation*}
\textbf{c}_1  < \textbf{c}_2 <  \textbf{c}_3.
\end{equation*}
We can now evaluate the formula from \cref{cor: geometric counting formula for inner product} to compute the entries of the matrix $\Psi$. For example, we have
\begin{equation}\label{eq: computation of some example of Psi}
\begin{aligned}
(E_{\textbf{i}, \textbf{a}_{\textbf{c}_2}}, E_{\textbf{i}, \textbf{a}_{\textbf{c}_2}} )&= v^{-\dim R_{\textbf{i}, \textbf{a}_{\textbf{c}_2}} - \dim R_{\textbf{i}, \textbf{a}_{\textbf{c}_2} }}  \tfrac{\prod_{i}(1-v^{-2})^{-\textbf{v}_i}}{[\textbf{a}_{\textbf{c}_2}]! [\textbf{a}_{\textbf{c}_2}]!}  \sum_{w \in S_2 \times S_2 } v^{2(l(w) + \dim R_{\textbf{i}, \textbf{a}_{\textbf{c}_2}} \cap {}^w R_{\textbf{i}, \textbf{a}_{\textbf{c}_2}})}\\
& = v^{-2-2} \tfrac{(1-v^{-2})^{-4}}{(v+v^{-1})(v+v^{-1}) } (v^{2(0+2)} + v^{2(1 + 2)} + v^{2(1+0)} + v^{2(2+0)}) \\
& = (1-v^{-2})^{-4} \\
& = \tfrac{1 + 2v^{-2}+ v^{-4}}{(1-v^{-2})^{2} (1-v^{-4})^{2}}.
\end{aligned}
\end{equation}
More generally, we have
\begin{equation*}
\Psi =  \tfrac{1}{(1-v^{-2})^2 (1-v^{-4})^2} \begin{spm}
1 & v^{-1} + v^{-3} & v^{-4} \\
v^{-1} + v^{-3} &  1 + 2v^{-2}+ v^{-4} & v^{-1} + v^{-3} \\
v^{-4} & v^{-1} + v^{-3} & 1
\end{spm}.
\end{equation*}
The matrices $L$ and $D$ in the corresponding LDLT decomposition $\Psi = LDL^T$ are
\begin{align*}
L &= \begin{spm}
1 & 0 & 0 \\ v^{-1} + v^{-3} & 1 & 0 \\ v^{-4} & v^{-1} & 1\\ 
\end{spm} \\
D &=  \tfrac{1}{(1-v^{-2})^{2} (1-v^{-4})^{2}} \begin{spm}
1 & 0 & 0 \\ 0 & 1+v^{-2}-v^{-4}-v^{-6} & 0 \\ 0 & 0 & 1+v^{-2}-v^{-4}-v^{-6} \end{spm} .
\end{align*} 
Note that $L$ already has coefficients in $v^{-1}\mathbb{Z}[v^{-1}]$ below the diagonal. Hence the change of basis matrix between $\Bstdv$ and $\Bcanv$ is
\begin{equation*}
P = L = \begin{spm}
1 & 0 & 0 \\ v^{-1} + v^{-3} & 1 & 0 \\ v^{-4} & v^{-1} & 1\\ 
\end{spm}.
\end{equation*}
The matrix $Q$ is just the identity matrix which corresponds to the fact that
\begin{equation*}
\Bcanv = \Bmonv = \{E_2^{(2)}E_1^{(2)}, E_1E_2^{(2)}E_1, E_1^{(2)} E_2^{(2)} \}.
\end{equation*}
\begin{remark}
The matrix $\Psi$ can also be computed recursively from \cite[1.2.3]{lusztig2010introduction} or using the closed formula in \cite{reineke2001monomials} derived from this recursion. To illustrate the difference to our formula, we briefly explain how to compute $\Psi$ with this formula. For any $a \in \mathbb{N}$ let
\begin{equation*}
(\{a\})^{-1}= \prod_{s=1}^a (1-v^{-2s})^{-1}.
\end{equation*}
The formula from \cite[Theorem 2.2]{reineke2001monomials} is then
\begin{equation}\label{eq: Reineckes formula}
(E_{\textbf{i}, \textbf{a}}, E_{\textbf{i}', \textbf{a}'} ) = \sum_{z \in \mathcal{Z}^{\textbf{i}, \textbf{a}}_{\textbf{i}', \textbf{a}'}} v^{Q(z)} \prod_{m,p} (\{ z_{mp} \})^{-1}
\end{equation}
where $\mathcal{Z}^{\textbf{i}, \textbf{a}}_{\textbf{i}', \textbf{a}'}$ is a certain finite set of $|\textbf{i}'| \times |\textbf{i}|$ matrices with entries in $\mathbb{N}$ and $Q(z)$ is a certain integer which can be computed from the entries of $z$. For example, one can check that
\begin{equation*}
\mathcal{Z}^{\textbf{i}, \textbf{a}_{\textbf{c}_2}}_{\textbf{i}, \textbf{a}_{\textbf{c}_2}} = \left\{ \begin{spm} 1 & 0 & 0 \\ 0 & 2 & 0 \\ 0 & 0 & 1 \end{spm}, \begin{spm} 0 & 0 & 1 \\ 0 & 2 & 0 \\ 1 & 0 & 0 \end{spm} \right\}
\end{equation*}
and
\begin{equation*}
Q\left( \begin{spm} 1 & 0 & 0 \\ 0 & 2 & 0 \\ 0 & 0 & 1 \end{spm} \right) = 0, \quad Q\left( \begin{spm} 0 & 0 & 1 \\ 0 & 2 & 0 \\ 1 & 0 & 0 \end{spm} \right) = -2.
\end{equation*}
Hence, we get
\begin{align*}
(E_{\textbf{i}, \textbf{a}_{\textbf{c}_2}}, E_{\textbf{i}, \textbf{a}_{\textbf{c}_2}} )&\overset{\eqref{eq: Reineckes formula}}{=} v^{0}(\{1\})^{-2} (\{2\})^{-1} + v^{-2}(\{1\})^{-2} (\{2\})^{-1} \\
& = \tfrac{1}{(1-v^{-2})^3(1-v^{-4})} + v^{-2}\tfrac{1}{(1-v^{-2})^3(1-v^{-4})} \\
&= (1-v^{-2})^{-4}
\end{align*}
which is the same result as in \eqref{eq: computation of some example of Psi}.
\end{remark}

\subsection{Example: \texorpdfstring{$\textbf{v} = (1,2,1)$}{v = (1,2,1)}}\label{section: composition multiplicities for 121}
Let $n=3$ and $\textbf{v} = (1,2,1)$. The order on the set of positive roots is
\begin{equation*}
\Phi^+  = \{ \textbf{e}_3 < \textbf{e}_2 + \textbf{e}_3 < \textbf{e}_2 < \textbf{e}_1 + \textbf{e}_2 + \textbf{e}_3 < \textbf{e}_1 + \textbf{e}_2 < \textbf{e}_1   \}.
\end{equation*}
The Kostant partitions of dimension $\textbf{v}$ are:
\begin{equation}\label{eq: kostant partitions of type 1 2 1}
\scriptstyle (1,0,2,0,0,1) , \quad  (1,0,1,0,1,0) , \quad (0,1,1,0,0,1)  , \quad (0,1,0,0,1,0) , \quad (0,0,1,1,0,0).
\end{equation}
We have
\begin{equation*}
\textbf{i} = (3,2,3,1,2,3)
\end{equation*}
and the corresponding sequences $a_{\textbf{c}}$ are
\begin{equation*}
\scriptstyle (1, 2, 0, 1,0,0), \quad  (1, 1, 0, 1,1,0),  \quad (0, 2, 1, 1,0,0), \quad (0, 1, 1, 1,1,0), \quad  (0, 1, 0, 1,1,1).
\end{equation*}
Let $B \subset \GL_{\textbf{v}}$ be the Borel consisting of all lower triangular matrices. Then the corresponding parabolic subgroups $P_{\textbf{i}, \textbf{a}_{\textbf{c}}}$ are
\begin{equation*}
\GL_{\textbf{v}} ,\quad  B, \quad  \GL_{\textbf{v}}, \quad  B,  \quad B.
\end{equation*}
The corresponding spaces $R_{\textbf{i}, \textbf{a}_{\textbf{c}}} \subset R_{\textbf{v}} = (\begin{spm} * \\ * \end{spm}, \begin{spm} * & * \end{spm} )$ are
\begin{equation}\label{eq: spaces R i a for v eq 121}
(\begin{spm} 0 \\ 0 \end{spm}, \begin{spm} 0 & 0 \end{spm} ), \quad (\begin{spm} 0 \\ * \end{spm} , \begin{spm} 0 & 0 \end{spm}), \quad  (\begin{spm} 0 \\ 0 \end{spm}, \begin{spm} * & * \end{spm} ), \quad   (\begin{spm} 0 \\ * \end{spm} , \begin{spm} * & 0 \end{spm} ), \quad   (\begin{spm} * \\ 0 \end{spm} ,\begin{spm} * & * \end{spm}).
\end{equation}
It follows from \eqref{eq: orbit dimension formula} that the corresponding orbit dimensions are
\begin{equation*}
0, \quad 2, \quad 2, \quad 3, \quad 4.
\end{equation*}
Thus the Kostant partitions in \eqref{eq: kostant partitions of type 1 2 1} are already ordered with respect to the partial order. The entries of the matrix $\Psi$ can now be computed explicitly. For example, taking
\begin{equation*}
\textbf{c}_3 = (0,1,1,0,0,1), \quad \textbf{c}_5 = (0,0,1,1,0,0)
\end{equation*}
we get
\begin{align*}
\Psi_{\textbf{c}_3, \textbf{c}_5} &= v^{- \dim R_{\textbf{i}, \textbf{a}_3} - \dim R_{\textbf{i}, \textbf{a}_5} }  \tfrac{\prod_i (1-v^{-2})^{-\textbf{v}_i}}{[\textbf{a}_{\textbf{c}_3}]! [\textbf{a}_{\textbf{c}_5}]!} (v^{2\dim R_{\textbf{i}, \textbf{a}_{\textbf{c}_3}} \cap R_{\textbf{i}, \textbf{a}_{\textbf{c}_5}}} + v^{2+ 2 \dim R_{\textbf{i}, \textbf{a}_{\textbf{c}_3}} \cap {}^s R_{\textbf{i}, \textbf{a}_{\textbf{c}_5}}} ) \\
&= v^{-5} \tfrac{(1-v^{-2})^{-4}}{[2]! \cdot 1} (v^4 + v^6)\\
&= (1-v^{-2})^{-4} \\
&= \tfrac{1+v^{-2}}{(1-v^{-2})^{-3}(1-v^{-4})}.
\end{align*}
More generally, we have
\begin{equation}\label{eq: psi in the case 121}
\Psi = \tfrac{1}{(1-v^{-2})^{3}(1-v^{-4})} \begin{spm}
1 & 1+v^{-2} & v^{-2} & v^{-1}+v^{-3} & v^{-2}+v^{-4} \\
1+v^{-2} & 2+2v^{-2} & v^{-2}+v^{-4} & 2v^{-1}+2v^{-3} &2v^{-2}+2v^{-4} \\
v^{-2} & v^{-2}+v^{-4} & 1 & v^{-1} +v^{-3} & 1+v^{-2}\\
v^{-1}+v^{-3} & 2v^{-1} + 2v^{-3} & v^{-1}+ v^{-3} & 1+ 2v^{-2} + v^{-4} & 2v^{-1} + 2v^{-3} \\
v^{-2} + v^{-4} & 2v^{-2} + 2v^{-4} & 1+v^{-2} & 2v^{-1} + 2v^{-3} & 2+2v^{-2}
\end{spm}.
\end{equation}
The decomposition $\Psi = LDL^T$ is given by
\begin{align*}
L &=\begin{spm}
1 & 0 & 0 & 0 & 0 \\
1+v^{-2} & 1 & 0 & 0 & 0 \\
v^{-2} & 0 & 1 & 0 & 0 \\
v^{-1}+ v^{-3} & v^{-1} & v^{-1} & 1 & 0 \\
v^{-2} + v^{-4} & v^{-2} & 1 + v^{-2}& v^{-1} & 1 
\end{spm} \\
D &=  \tfrac{1}{(1-v^{-2})^{3}(1-v^{-4})}
\begin{spm}
1  \\
 & 1-v^{-4}  \\
 &  & 1-v^{-4}  \\
 &  &  & 1-v^{2}-v^{-4}+v^{-6} \\
 &  &  &  & 1-v^{-2}-v^{-4}+v^{-6} 
\end{spm}.
\end{align*}
The matrix $L$ can be decomposed further as
\begin{equation}\label{eq: matrix p and q in the v is 121 example}
L = QP =  \begin{spm}
1 & 0 & 0 & 0 & 0\\
1 & 1 & 0 & 0 & 0\\
0 & 0 & 1 & 0 & 0\\
0 & 0 & 0 & 1 & 0\\
0 & 0 & 1 & 0 & 1
\end{spm}
\begin{spm}
1 & 0 & 0 & 0 & 0 \\
v^{-2} & 1 & 0 & 0 & 0 \\
v^{-2} & 0 & 1 & 0 & 0 \\
v^{-1} + v^{-3} & v^{-1} & v^{-1} & 1 & 0 \\
v^{-4} & v^{-2} & v^{-2}& v^{-1} & 1 
\end{spm}.
\end{equation}
\begin{remark}
The matrix $P$ is also computed in the context of (graded) affine Hecke algebras in \cite[4.2]{ciubotaru2008multiplicity} using the slightly different normalisation
\begin{equation*}
\tilde{P} = ( v^{\dim \mathcal{O}_{\textbf{c}'}-\dim \mathcal{O}_{\textbf{c}}}p_{\textbf{c}',\textbf{c}}(v) )_{\textbf{c},\textbf{c}' \in \KP(\textbf{v})}  = \begin{spm}
1 & 1 & 1 & 1+ v^2 & 1 \\
0 & 1 & 0 & 1 & 1 \\
0 & 0 & 1 &  1 & 1 \\
0 & 0 & 0 & 1 & 1 \\
0 & 0 & 0 & 0 & 1
\end{spm}.
\end{equation*}
The relation between the quantum groups setting and the affine Hecke algebra setting will be discussed in more detail in the next section.
\end{remark}

\section{Applications to the representation theory of affine Hecke algebras}
\subsection{Affine Hecke algebras and quivers}\label{section: hecke algebras and quivers}
Let $\mathcal{H}^{\aff} = \mathcal{H}^{\aff}(\GL_n)$ be the affine Hecke algebra of $\GL_n$. This is the $\mathbb{C}[q,q^{-1}]$-algebra with generators
\begin{equation*}
T_1 , ..., T_{n-1}, \quad  X_1^{\pm 1}, ..., X_n^{\pm 1}
\end{equation*}
and relations
\begin{itemize}
\item[]
\item[] \makebox[7cm]{$(T_i+1)(T_i-q) = 0 $ \hfill}  $i = 1 ,...,n-1$;
\item[] \makebox[7cm]{$ T_iT_j = T_j T_i $ \hfill}  $| i - j | > 1$;
\item[] \makebox[7cm]{$T_i T_{i+1} T_i = T_{i+1}T_i T_{i+1} $ \hfill} $i = 1, ..., n-2$;
\item[]
\item[] \makebox[7cm]{$ X_i X_i^{-1} = 1 = X_i^{-1} X_i $ \hfill}  $i= 1 , ..., n$;
\item[] \makebox[7cm]{$X_i X_j  =  X_j X_i $ \hfill} for all $i,j$;
\item[]
\item[] \makebox[7cm]{$T_i X_j  =  X_j T_i $ \hfill} $|i-j| > 1$;
\item[] \makebox[7cm]{$T_iX_i T_i = qX_{i+1} $ \hfill} $i = 1, ..., n-1$.
\item[]
\end{itemize}
The center of the affine Hecke algebra can be computed as
\begin{equation*}
Z(\mathcal{H}^{\aff}) \cong \mathbb{C}[q,q^{-1}][X_1^{\pm 1}, ..., X_n^{\pm 1}]^{S_n} \cong  \mathcal{O}(\GL_n \times \mathbb{C}^{\times})^{\GL_n}  .
\end{equation*}
As a consequence, the central characters of $\mathcal{H}^{\aff}$ are parametrized by semisimple conjugacy classes in $ \GL_n \times \mathbb{C}^{\times}$. For any $(s,q_0) \in \GL_n \times \mathbb{C}^{\times}$ semisimple, we denote the corresponding central character by
\begin{equation*}
\chi_{(s,q_0)} : Z( \mathcal{H}^{\aff}) \rightarrow \mathbb{C}.
\end{equation*}
Any simple $\mathcal{H}^{\aff}$-module admits a central character by (the countable dimension version of) Schur's lemma, i.e. there is a unique central character $\chi_{(s,q_0)}$ such that $\ker(\chi_{(s,q_0)})$ annihilates the module. Hence, there is a decomposition
\begin{equation*}
\Irr(\mathcal{H}^{\aff}\modu) = \bigsqcup_{(s,q_0)/\sim } \Irr_{(s,q_0)}(\mathcal{H}^{\aff} \modu) 
\end{equation*}
where $\Irr_{(s,q_0)}(\mathcal{H}^{\aff} \modu) $ is the set of all simple $\mathcal{H}^{\aff}$-modules with central character $\chi_{(s,q_0)}$. When $q_0$ is not a root of unity, the set $\Irr_{(s,q_0)}(\mathcal{H}^{\aff}\modu)$ can be parametrized using geometric techniques. Pick a Borel subgroup $B \subset \GL_n$ with Lie algebra $\mathfrak{b} \subset \mathfrak{gl}_n$. Let $\mathcal{N} \subset \mathfrak{gl}_n$ be the nilpotent cone and $\mathfrak{n} := [\mathfrak{b}, \mathfrak{b}]$. Let
\begin{equation*}
\tilde{\mathcal{N}} := \GL_n \times^B \mathfrak{n}
\end{equation*}
and consider the Springer resolution
\begin{align*}
\mu : \tilde{\mathcal{N}} &\rightarrow \mathcal{N} \\
 (g,n) & \mapsto g\cdot n.
\end{align*}
This morphism is $\GL_n \times\mathbb{C}^{\times}$-equivariant where $t \in \mathbb{C}^{\times}$ acts by scaling with $t^{-1}$ on $\mathcal{N}$ (resp. $\mathfrak{n}$). Taking $(s,q_0)$-fixed points, we obtain a morphism
\begin{equation*}
\mu^{(s,q_0)} : \tilde{\mathcal{N}}^{(s,q_0)} \rightarrow \mathcal{N}^{(s,q_0)}
\end{equation*}
which is $\GL_n(s)$-equivariant. Explicitly, we have
\begin{equation*}
\mathcal{N}^{(s,q_0)} = \{ x \in \mathcal{N} \mid sxs^{-1} = q_0x \}.
\end{equation*}
For any $x \in \mathcal{N}^{(s,q_0)}$ let
\begin{equation*}
\mathcal{O}_x := \GL_n(s) \cdot x \subset \mathcal{N}^{(s,q_0)}
\end{equation*}
be the corresponding orbit and let
\begin{equation*}
M_x := H_*((\mu^{(s,q_0)})^{-1}(x))
\end{equation*}
be the Borel-Moore homology of the corresponding fiber. The following theorem is also known as the Deligne-Langlands correspondence for $\GL_n$.
\begin{theorem}\label{thm: deligne-langlands correspondence}\cite{bernstein1977induced,zelevinsky1980induced,kazhdan1987proof,chriss2009representation}
Assume that $q_0 \in \mathbb{C}^{\times}$ is not a root of unity.
\begin{enumerate}
\item There is a canonical $\mathcal{H}^{\aff}$-module structure on $M_x$ with central character $\chi_{(s,q_0)}$.
\item If $\mathcal{O}_x = \mathcal{O}_{x'}$, then $M_x \cong M_{x'}$ as $\mathcal{H}^{\aff}$-modules.
\item The module $M_x$ has a unique simple quotient $L_x$. This induces a bijection
\begin{align*}
\{ \GL_n(s) \text{-orbits on } \mathcal{N}^{(s,q_0)} \} & \overset{1:1}{\leftrightarrow} \Irr_{(s,q_0)}( \mathcal{H}^{\aff} \modu )  \\
\mathcal{O}_x & \mapsto L_x.
\end{align*}
\item The composition multiplicities of the $M_x$ are given by the formula
\begin{equation*}
[M_x : L_y] = \sum_{k \in \mathbb{Z}} \dim H^k ( \iota_x^! \IC_{\mathcal{O}_y}) = \sum_{k \in \mathbb{Z}} \dim H^k ( \iota_x^* \IC_{\mathcal{O}_y})
\end{equation*}
where $\iota_x : \{x \} \hookrightarrow \mathcal{N}^{(s,q_0)}$ is the inclusion.
\end{enumerate}
\end{theorem}
It turns out that $\mathcal{N}^{(s,q_0)}$ is just a space of quiver representations in disguise. Assume for simplicity that $s$ is a diagonal matrix of the form
\begin{equation}\label{eq: definition of diag matrix with q powers}
s = s_{e, \textbf{v}} := e \cdot \begin{pmatrix}
\Id_{\textbf{v}_1} \\ & q_0 \Id_{\textbf{v}_2} \\ & & \ddots \\ & & &  q_0^{k-1} \Id_{\textbf{v}_k}
\end{pmatrix}
\end{equation}
for some integers $\textbf{v}_1, ..., \textbf{v}_m \in \mathbb{N}$ and $e, q_0 \in \mathbb{C}^{\times}$ where $q_0$ is not a root of unity (see \cref{remark: general semisimple element} for the case of a general semisimple $s$). Then we have
\begin{equation}\label{eq: centraliser of semisimple element in GLn}
\begin{aligned}
\GL_n(s) & = \begin{pmatrix}
\GL_{\textbf{v}_1} \\ & \GL_{\textbf{v}_2} \\ & & \ddots \\ & 
& &  \GL_{\textbf{v}_m} 
\end{pmatrix}
 \\
\mathcal{N}^{(s,q_0)} & =  {\tiny
\left( \begin{array}{ccccc}
0_{\textbf{v}_1 \times \textbf{v}_1} &&&& \\  \\
\Hom(\mathbb{C}^{\textbf{v}_1}, \mathbb{C}^{\textbf{v}_2}) & 0_{\textbf{v}_2 \times \textbf{v}_2} &&&  \\ \\
& \Hom(\mathbb{C}^{\textbf{v}_2}, \mathbb{C}^{\textbf{v}_3}) & 0_{\textbf{v}_3 \times \textbf{v}_3} && \\ \\
& &  \downdots  & \downdots &  \\ \\
&&& \Hom(\mathbb{C}^{\textbf{v}_{m-1}}, \mathbb{C}^{\textbf{v}_m}) & \ 0_{\textbf{v}_m \times \textbf{v}_m}
\end{array} \right) }.
\end{aligned}
\end{equation}
Thus, we can identify
\begin{align*}
\GL_n(s) & \cong \GL_{\textbf{v}} \\
 \mathcal{N}^{(s,q_0)} & \cong R_{\textbf{v}}(A_m)
\end{align*}
where $A_m$ is the quiver from \eqref{eq: type A quiver}. Hence, the sets of simple and standard modules with central character $\chi_{(s,q_0)}$ are indexed by Kostant partitions:
\begin{equation*}
\{ L_{\textbf{c}} \mid \textbf{c} \in \KP(\textbf{v}) \} , \quad \{ M_{\textbf{c}} \mid \textbf{c} \in \KP(\textbf{v}) \}.
\end{equation*}
Let $U^+$ be the positive part of the quantum group associated to $A_m$.
\begin{lemma}
We have $[M_{\textbf{c}}: L_{\textbf{c}'}] = p_{\textbf{c}', \textbf{c}}(1)$ where $P = (p_{\textbf{c}, \textbf{c}'}(v))_{\textbf{c} , \textbf{c}'  \in \KP(\textbf{v})}$ is the change of basis matrix $P$ between $\Bcanv$ and $\Bstdv$ in $U^+_{\textbf{v}}$.
\end{lemma}
\begin{proof}
By \cref{thm: deligne-langlands correspondence} and \cref{prop: geometric interpretation of base change between standard and canonical basis} we have
\begin{equation*}
[M_{\textbf{c}}: L_{\textbf{c}'}] = \sum_{k \in \mathbb{Z}} \dim H^k (\iota_x^* \IC_{\textbf{c}'}) = p_{\textbf{c}', \textbf{c}}(1)
\end{equation*}
where $\iota_x : \{x \} \hookrightarrow \mathcal{O}_{\textbf{c}}$ is any base point.
\end{proof}
\begin{corollary}
The multiplicities $[M_{\textbf{c}}: L_{\textbf{c}'}]$ can be computed with \cref{algorithm: new algo}.
\end{corollary}
\begin{example}\label{example: composition multiplicities for 121}
Let $\textbf{v} = (1,2,1)$ and consider a semisimple element of the form
\begin{equation*}
s= s_{e, \textbf{v}} = \Diag (e,  eq_0, eq_0,  eq_0^2  ) \in \GL_4.
\end{equation*}
Then
\begin{equation*}
\GL_4(s) = \GL_{(1,2,1)}, \quad \mathcal{N}^{(s,q_0)} = R_{(1,2,1)}
\end{equation*}
is the geometry considered in \cref{section: composition multiplicities for 121}. In particular, there are five orbits and hence $\mathcal{H}^{\aff}$ has five simple (resp. standard) modules with central character $\chi_{(s,q_0)}$
\begin{equation*}
L_1, ..., L_5, \quad M_1, ..., M_5.
\end{equation*}
Recall from \eqref{eq: matrix p and q in the v is 121 example} that the base change of basis matrix $P$ between $\Bstdv$ and $\Bcanv$ in $U^+_{\textbf{v}}$ is given by
\begin{equation*}\\
P = \begin{spm}
1 & 0 & 0 & 0 & 0 \\
v^{-2} & 1 & 0 & 0 & 0 \\
v^{-2} & 0 & 1 & 0 & 0 \\
v^{-1} + v^{-3} & v^{-1} & v^{-1} & 1 & 0 \\
v^{-4} & v^{-2} & v^{-2}& v^{-1} & 1 
\end{spm}
\end{equation*}
Hence, we get
\begin{equation*}
[M_i : L_j] = p_{j,i}(1) = \begin{cases}
0 & i >  j \\
0 & i = 2, j = 3\\
2 & i = 1, j=4\\
1 & \text{otherwise}.
\end{cases}
\end{equation*}
\end{example}
\begin{remark}\label{remark: general semisimple element}
For a general semisimple pair $(s,q_0) \in \GL_n \times \mathbb{C}^{\times}$ such that $q_0$ is not a root of unity, one can still show that there exist integers $m_1 , ..., m_k$ and sequences $\textbf{v}^1 ,..., \textbf{v}^k$ (where $\textbf{v}^i \in \mathbb{N}^{m_i}$) such that
\begin{equation}\label{eq: structure of centraliser in the general case}
\begin{aligned}
\GL_n(s) & \cong \GL_{\textbf{v}^1} \times ... \times \GL_{\textbf{v}^k}\\
\mathcal{N}^{(s,q_0)} & \cong R_{\textbf{v}^1}(A_{m_1}) \times ... \times R_{\textbf{v}^k}(A_{m_k}).
\end{aligned}
\end{equation}
In fact, we may conjugate $s$ into a diagonal matrix of the form
\begin{equation*}
s = \Diag (s_{e_1 , \textbf{v}^1}, s_{e_2, \textbf{v}^2} , ..., s_{e_k, \textbf{v}^k})
\end{equation*}
for some $e_i \in \mathbb{C}^{\times}$ and $\textbf{v}^i \in \mathbb{N}^{m_i}$ where $s_{e_i , \textbf{v}^i}$ is as in \eqref{eq: definition of diag matrix with q powers}. Moreover, we may assume that
\begin{equation*}
e_i/e_j \not\in \{ q_0^{m_i}, q_0^{m_i-1}, ..., q_0^{-m_j} \} \quad \text{for all } i \neq j
\end{equation*}
since otherwise $s_{e_i , \textbf{v}^i}$ and $s_{e_j ,\textbf{v}^j}$ can be combined into a single block. Then one can construct the isomorphisms in \eqref{eq: structure of centraliser in the general case} by the same argument as in \eqref{eq: centraliser of semisimple element in GLn}. When $q_0$ is a root of unity, the space $\mathcal{N}^{(s,q_0)}$ is related to representations of the cyclic quiver and affine Schubert varieties \cite[\S 11]{lusztig1990canonical}.
\end{remark}

\subsection{Dimensions of simple modules}\label{Section: dimensions of simples}
We now want to explain how to compute the dimensions of the simple $\mathcal{H}^{\aff}$-modules with central character $\chi_{(s,q_0)}$ for $s = s_{e, \textbf{v}}$ as in \eqref{eq: definition of diag matrix with q powers} (the case of a general semisimple element can be reduced to this case as in \cref{remark: general semisimple element}). Using the map
\begin{equation*}
\mu^{(s,q_0)} : \tilde{\mathcal{N}}^{(s,q_0)} \rightarrow \mathcal{N}^{(s,q_0)} \cong R_{\textbf{v}}(A_m).
\end{equation*}
one can define the sheaf
\begin{equation*}
\Sbf^{(s,q_0)} := (\mu^{(s,q_0)})_* \mathcal{C}_{\tilde{\mathcal{N}}^{(s,q_0)}} \in D^b_{\GL_{\textbf{v}}}(R_{\textbf{v}}(A_m)).
\end{equation*}
By the decomposition theorem, we have
\begin{equation*}
\Sbf^{(s,q_0)} = \sum_{\textbf{c} \in \KP(\textbf{v})} f'_{\textbf{c}}(v) \cdot \IC_{\textbf{c}}
\end{equation*}
for some $f'_{\textbf{c}}(v) \in \mathbb{Z}[v,v^{-1}]$.
\begin{lemma}\label{lem: dim of simple is multiplicities of perverse sheaves}\cite[8.6.12]{chriss2009representation}
We have $\dim L_{\textbf{c}} = f'_{\textbf{c}}(1)$.
\end{lemma}
Thus to compute the dimensions of simple modules, it is sufficient to compute the polynomials $f'_{\textbf{c}}(v)$. For this we need to analyse the geometry of $\tilde{\mathcal{N}}^{(s,q_0)}$ in a bit more detail. Let $B \subset \GL_n$ be a Borel subgroup containing $s$, $T \subset B$ a maximal torus containing $s$ and $W=S_n$ the corresponding Weyl group. Then $B(s)$ is a Borel subgroup of $\GL_n(s)$ \cite[8.8.7]{chriss2009representation}.
\begin{lemma}
The connected components of $\tilde{\mathcal{N}}^{(s,q_0)}$ are of the form
\begin{equation*}
\tilde{\mathcal{N}}^{(s,q_0), w} := G(s)\times^{B(s)} (  {}^{w} \mathfrak{n} \cap \mathcal{N}^{(s,q_0)})
\end{equation*}
where $w$ runs through the set ${}^{\textbf{v}} S_n$ of minimal (right) coset representatives in $S_{\textbf{v}} \backslash S_n$.
\end{lemma}
\begin{proof}
Let $\mathcal{B} = \GL_n/B$ be the flag variety of $\GL_n$ and $\mathcal{B}^s$ the variety of $s$-fixed points, i.e. the variety of all Borel subgroups containing $s$. Let
\begin{equation*}
\pi^{(s,q_0)} : \tilde{\mathcal{N}}^{(s,q_0)} \rightarrow \mathcal{B}^s
\end{equation*}
be the projection. By \cite[8.8.7]{chriss2009representation} the connected components of $\mathcal{B}^s$ are all isomorphic to the flag variety $\GL_n(s)/B(s)$ of $\GL_n(s)$. Thus each connected component $C$ of $\mathcal{B}^s$ contains a unique Borel subgroup $B'$ whose stabilizer in $G(s)$ is $B(s)$, i.e. $B'(s) = B(s)$. By \cite[8.8.9]{chriss2009representation} the connected components of $\tilde{\mathcal{N}}^{(s,q_0)}$ are of the form
\begin{equation*}
(\pi^{(s,q_0)})^{-1}( C) \cong \GL_n(s) \times^{B(s)} (\mathfrak{n}' \cap \mathcal{N}^{(s,q_0)}).
\end{equation*}
Thus, to complete the proof, it remains to show that
\begin{equation}\label{eq description of borels with B'(s) = B(s)}
\{ B' \in  \mathcal{B}^s \mid B'(s) = B(s) \} = \{{}^wB \mid w \in {}^{\textbf{v}} S_n \}.
\end{equation}
To see this, note that any $B'$ with $B'(s) = B(s)$ contains the maximal torus $T \subset B(s)$ and thus $B'= {}^wB$ for some $w \in S_n$. Let $\Phi$ be the root system of $\GL_n$ and let
\begin{equation*}
\Phi_{\textbf{v}} := \{ \alpha \in \Phi \mid \alpha(s) =1 \}.
\end{equation*}
This a subrootsystem of $\Phi$ whose Weyl group is the parabolic subgroup $S_{\textbf{v}} \subset S_n$. The set of non-zero roots of $B'(s) = {}^wB \cap \GL_n(s)$ is
\begin{equation*}
w(\Phi^+) \cap \Phi_{\textbf{v}} .
\end{equation*}
Thus, we have $B'(s) = B(s)$ if and only if 
\begin{align*}
&w(\Phi^+) \cap \Phi_{\textbf{v}}   = \Phi^+ \cap \Phi_{\textbf{v}} \\
\Leftrightarrow & \Phi^+ \cap \Phi_{\textbf{v}} \subset w(\Phi^+)\\
\Leftrightarrow& w^{-1}(\Phi^+ \cap \Phi_{\textbf{v}}) \subset \Phi^+.
\end{align*}

This is equivalent to $w^{-1}$ being a minimal (left) coset representative in $S_n/S_{\textbf{v}}$ or to $w$ being a minimal (right) coset representative in $S_{\textbf{v}} \backslash S_n$. Hence, we have verified \eqref{eq description of borels with B'(s) = B(s)} completing the proof.
\end{proof}
For any $w \in {}^{\textbf{v}} S_n$ let
\begin{equation*}
\Sbf^w := (\mu^{(s,q_0)})_* \mathcal{C}_{\tilde{\mathcal{N}}^{(s,q_0),w}} \in D^b_{\GL_{\textbf{v}}}(R_{\textbf{v}}(A_m))
\end{equation*}
and write
\begin{equation}\label{eq: def of fw}
\Sbf^{w} = \sum_{\textbf{c} \in \KP(\textbf{v})} f^w_{\textbf{c}}(v) \cdot \IC_{\textbf{c}}.
\end{equation}
Then
\begin{equation}\label{eq: dec of Sbf and f into w cosets}
\Sbf^{(s,q_0)} = \bigoplus_{w \in {}^{\textbf{v}}S_n} \Sbf^w  \quad \text{and} \quad f'_{\textbf{c}}(v) = \sum_{w \in {}^{\textbf{v}}S_n} f^w_{\textbf{c}}(v).
\end{equation}
We also define polynomials
\begin{equation}\label{eq: def of hw}
h^w_{\textbf{c}}(v) = \sum_{k \in \mathbb{Z}} \dim \Hom^k( \Sbf_{\textbf{i}, \textbf{a}_{\textbf{c}}}, \Sbf^w) v^{-k}
\end{equation}
for any $\textbf{c} \in \KP(\textbf{v})$. To get slightly nicer formulas later on, we also define renormalized polynomials:
\begin{equation}\label{eq: renormalised f and h}
\begin{aligned}
f_{\textbf{c}}(v) &:= \sum_{w \in {}^{\textbf{v}} S_n}  v^{2l(w) + \dim {}^{w} \mathfrak{n} \cap \mathcal{N}^{(s,q_0)}} f^w_{\textbf{c}}(v)\\
h_{\textbf{c}}(v) &:= \sum_{w \in {}^{\textbf{v}} S_n}  v^{2l(w) + \dim {}^{w} \mathfrak{n} \cap \mathcal{N}^{(s,q_0)}} h^w_{\textbf{c}}(v).
\end{aligned}
\end{equation}
Note that
\begin{equation*}
\dim L_{\textbf{c}} \overset{\ref{lem: dim of simple is multiplicities of perverse sheaves}}{=} f'_{\textbf{c}}(1)  \overset{\eqref{eq: dec of Sbf and f into w cosets}}{=}  \sum_{w \in {}^{\textbf{v}}S_n}  f^w_{\textbf{c}}(1) \overset{\eqref{eq: renormalised f and h}}{=} f_{\textbf{c}}(1)
\end{equation*}
so for the purpose of determining the dimensions of simple modules, it is sufficient to compute the renormalized polynomials $f_{\textbf{c}}(v)$. For any $w \in S_n$ we can consider the intersection $R_{\textbf{i}, \textbf{a}_{\textbf{c}}} \cap {}^w \mathfrak{n}$ by viewing $R_{\textbf{i}, \textbf{a}_{\textbf{c}}}$ as a subspace of $\mathfrak{g}$ via
\begin{equation*}
R_{\textbf{i}, \textbf{a}_{\textbf{c}}} \subset R_{\textbf{v}} \cong \mathcal{N}^{(s,q_0)} \subset \mathfrak{g}.
\end{equation*}
Here we assume that $R_{\textbf{i}, \textbf{a}_{\textbf{c}}}$ is defined with respect to the Borel  subgroup $B(s) \subset \GL_n(s) \cong\GL_{\textbf{v}}$.
\begin{lemma}\label{lem: formula for h}
We have $h_{\textbf{c}}(v) = v^{-\dim R_{\textbf{i}, \textbf{a}_{\textbf{c}}}} \frac{(1-v^{-2})^{-n}}{[\textbf{a}_{\textbf{c}}]!} \sum_{w \in S_n} v^{2(l(w) + \dim R_{\textbf{i}, \textbf{a}_{\textbf{c}}} \cap {}^{w} \mathfrak{n}) }$.
\end{lemma}
\begin{proof}
The parabolic $P_{\textbf{i}, \textbf{a}_{\textbf{c}}}$ contains the Borel $B(s)$. Hence, for any $w \in {}^{\textbf{v}}S_n$ we can apply \cref{theorem: Ext formula} with $G = \GL_n(s)$, $V = \mathcal{N}^{(s,q_0)} \cong R_{\textbf{v}}$,  $V_1 = R_{\textbf{i}, \textbf{a}_{\textbf{c}}}$ , $P_1 = P_{\textbf{i}, \textbf{a}_{\textbf{c}}}$, $V_2 =  {}^{w} \mathfrak{n} \cap \mathcal{N}^{(s,q_0)}$ and $P_2 = B(s)$ to obtain
\begin{equation}\label{eq: computation of hw}
\begin{aligned}
h^w_{\textbf{c}}(v) & = v^{- \dim \GL_n(s)/B(s) - \dim R_{\textbf{i}, \textbf{a}_{\textbf{c}}} - \dim {}^{w} \mathfrak{n} \cap \mathcal{N}^{(s,q_0)}} \tfrac{\chi^{\GL_n(s)}_v(\GL_n(s)/ B(s))}{\chi_v(P_{\textbf{i}, \textbf{a}_{\textbf{c}}}/B(s))} \\
& \quad \cdot \sum_{y \in S_{\textbf{v}}} v^{2(l(y) + \dim R_{\textbf{i}, \textbf{a}_{\textbf{c}}} \cap {}^y( {}^{w} \mathfrak{n} \cap \mathcal{N}^{(s,q_0)})} \\
& = v^{- \dim R_{\textbf{i}, \textbf{a}_{\textbf{c}}} - \dim {}^{w} \mathfrak{n} \cap \mathcal{N}^{(s,q_0)} } \tfrac{(1-v^{-2})^{-n}}{[\textbf{a}_{\textbf{c}}]!} \sum_{y \in S_{\textbf{v}}} v^{2(l(y) + \dim R_{\textbf{i}, \textbf{a}_{\textbf{c}}} \cap {}^{yw} \mathfrak{n})}.
\end{aligned}
\end{equation}
Hence, we get
\begin{align*}
h_{\textbf{c}}(v) & \overset{\eqref{eq: renormalised f and h}}{=} \sum_{w \in {}^{\textbf{v}} S_n}  v^{2l(w) + \dim {}^{w} \mathfrak{n} \cap \mathcal{N}^{(s,q_0)}} h^w_{\textbf{c}}(v) \\
& \overset{\eqref{eq: computation of hw}}{=}  \sum_{w \in {}^{\textbf{v}} S_n}  v^{ 2l(w) - \dim R_{\textbf{i}, \textbf{a}_{\textbf{c}}}} \tfrac{(1-v^{-2})^{-n}}{[\textbf{a}_{\textbf{c}}]!} \sum_{y \in S_{\textbf{v}}} v^{2(l(y) + \dim R_{\textbf{i}, \textbf{a}_{\textbf{c}}} \cap {}^{yw} \mathfrak{n})} \\
& =  \sum_{w \in {}^{\textbf{v}} S_n}  v^{ - \dim R_{\textbf{i}, \textbf{a}_{\textbf{c}}}} \tfrac{(1-v^{-2})^{-n}}{[\textbf{a}_{\textbf{c}}]!} \sum_{y \in S_{\textbf{v}}} v^{2(l(yw) + \dim R_{\textbf{i}, \textbf{a}_{\textbf{c}}} \cap {}^{yw} \mathfrak{n})} \\
& =  v^{ - \dim R_{\textbf{i}, \textbf{a}_{\textbf{c}}}} \tfrac{(1-v^{-2})^{-n}}{[\textbf{a}_{\textbf{c}}]!} \sum_{w \in S_n} v^{2(l(w) + \dim R_{\textbf{i}, \textbf{a}_{\textbf{c}}} \cap {}^{w} \mathfrak{n} )}.
\end{align*}
\end{proof}
Consider the column vectors
\begin{align*}
F & := (f_{\textbf{c}}(v))_{\textbf{c} \in \KP(\textbf{v})} \\
H & := (h_{\textbf{c}}(v))_{\textbf{c} \in \KP(\textbf{v})}.
\end{align*}
Recall that $\Psi$ is the matrix representing the standard bilinear form on $U^+_{\textbf{v}}$ with respect to $\Bmonv$ and $Q$ is the change of basis matrix between $\Bcanv$ and $\Bmonv$.
\begin{lemma}
We have $H = \Psi Q^{-T} F$.
\end{lemma}
\begin{proof}
By \cref{thm: geometric categorification of qua group} and \cref{prop: i a springer sheaf identifies with monomial} we have
\begin{equation}\label{eq: matrix Q encodes decomposition of the Sbf i a}
\Sbf_{\textbf{i}, \textbf{a}_{\textbf{c}}} = \sum_{\textbf{c}' \in \KP(\textbf{v})} q_{\textbf{c}', \textbf{c}}(v) \cdot \IC_{\textbf{c}'}.
\end{equation}
Hence, we get
\begin{align*}
h^w_{\textbf{c}}(v)   \overset{\eqref{eq: def of hw}}&{=} \sum_{k \in \mathbb{Z}} \dim \Hom^k( \Sbf_{\textbf{i}, \textbf{a}_{\textbf{c}}}, \Sbf^w) v^{-k} \\ 
\overset{\eqref{eq: def of fw}, \eqref{eq: matrix Q encodes decomposition of the Sbf i a}}&{=}  \sum_{k \in \mathbb{Z}}    \sum_{\textbf{c}',\textbf{c}'' }  \overline{q}_{\textbf{c}, \textbf{c}'}(v) f^w_{\textbf{c}''}(v)  \dim \Hom^k( \IC_{\textbf{c}'}, \IC_{\textbf{c}''}) v^{-k} \\
\overset{\eqref{eq: matrix Q is bar invariant}}&{=}  \sum_{k \in \mathbb{Z}}    \sum_{\textbf{c}',\textbf{c}''}  q_{\textbf{c}, \textbf{c}'}(v) f^w_{\textbf{c}''}(v) \dim \Hom^k(  \IC_{\textbf{c}'}, \IC_{\textbf{c}''})  v^{-k} \\
\overset{\ref{thm: geometric categorification of qua group}}&{=}    \sum_{\textbf{c}',\textbf{c}'' }  q_{\textbf{c}, \textbf{c}'}(v) \omega_{\textbf{c}', \textbf{c}''}(v) f^w_{\textbf{c}''}(v)
\end{align*}
where the $ \omega_{\textbf{c}', \textbf{c}''}(v)$ are the entries of the matrix $\Omega$ representing the standard bilinear form on $U^+_{\textbf{v}}$ with respect to $\Bcanv$. Thus, we get
\begin{align*}
h_{\textbf{c}}(v) & \overset{\eqref{eq: renormalised f and h}}{=} \sum_{w \in {}^{\textbf{v}} S_n}  v^{2l(w) + \dim {}^{w} \mathfrak{n} \cap \mathcal{N}^{(s,q_0)}} h^w_{\textbf{c}}(v) \\
& = \sum_{w \in {}^{\textbf{v}} S_n}  v^{2l(w) + \dim {}^{w} \mathfrak{n} \cap \mathcal{N}^{(s,q_0)}} \sum_{\textbf{c}',\textbf{c}'' \in \KP(\textbf{v}) }  q_{\textbf{c}, \textbf{c}'}(v) \omega_{\textbf{c}', \textbf{c}''}(v) f^w_{\textbf{c}''}(v) \\
& \overset{\eqref{eq: renormalised f and h}}{=}  \sum_{\textbf{c}',\textbf{c}'' \in \KP(\textbf{v})}  q_{\textbf{c}, \textbf{c}'}(v) \omega_{\textbf{c}', \textbf{c}''}(v) f_{\textbf{c}''}(v) .
\end{align*}
As a matrix equation, this is precisely $H =Q \Omega F$. Hence, we get
\begin{equation*}
H = Q\Omega Q^T Q^{-T} F \overset{\eqref{eq: Psi is Q omega Q t}}{=}  \Psi Q^{-T}F.
\end{equation*}
\end{proof}
Recall that we can compute $Q$ and $\Psi$ using \cref{algorithm: new algo}. The vector $H$ can be computed explicitly using \cref{lem: formula for h}. Hence, using the lemma above, we can compute $F$ from $Q, \Psi$ and $H$. Thus, we arrive at the following algorithm to compute the dimensions of simple $\mathcal{H}^{\aff}$-modules with central character $\chi_{(s,q_0)}$ where $s = s_{e, \textbf{v}}$.
\begin{algorithm}\label{algorithm: dimensions of simple modules}
Input: $\textbf{v}$, Output: $\dim L_{\textbf{c}}$ for all $\textbf{c} \in \KP(\textbf{v})$;\\
Step 1: Compute the matrices $\Psi$ and $Q$ using \cref{algorithm: new algo};\\
Step 2: Compute the vector $H$ using \cref{lem: formula for h};\\
Step 3: Compute the vector $F = Q^T \Psi^{-1}H$.\\
Step 4: Compute $\dim L_{\textbf{c}} = f_\textbf{c}(1)$.
\end{algorithm}
\begin{example}\label{example: dimensions of simple for v equals 121}
Let $\textbf{v} = (1,2,1)$ and $s= s_{e, \textbf{v}}$. We compute the dimensions of the simple $\mathcal{H}^{\aff}(\GL_4)$-modules with central character $\chi_{(s,q_0)}$. This continues the examples from \cref{section: composition multiplicities for 121} and \cref{example: composition multiplicities for 121}. We use the Borel $B \subset \GL_4$ consisting of lower triangular matrices. In \cref{table: intersections}, we list for each $w \in S_4$ (in cycle notation) the length of the element $w$ and the subspace
\begin{equation*}
\mathcal{N}^{(s,q_0)} \cap {}^w \mathfrak{n} =  R_{\textbf{v}} \cap {}^w \mathfrak{n}\subset  R_{\textbf{v}} = \left( \begin{spm} * \\ * \end{spm}, \begin{spm} * & * \end{spm} \right).
\end{equation*}
\begin{table}
\caption{\label{table: intersections}}
\begin{equation*}
\begin{array}{c|c|c}
w & l(w) &   R_{\textbf{v}} \cap {}^w \mathfrak{n} \\ \hline \hline
\begin{matrix} e \\ (23) \end{matrix} & \begin{matrix} 0 \\ 1 \end{matrix} & \left( \begin{spm} * \\ * \end{spm}, \begin{spm} * & * \end{spm} \right) \\ \hline
(12) & 1 &  \left( \begin{spm} 0 \\ * \end{spm}, \begin{spm} * & * \end{spm} \right)\\ \hline
(34) & 1 & \left( \begin{spm} * \\ * \end{spm}, \begin{spm} * & 0 \end{spm} \right) \\ \hline
(132) & 2 & \left( \begin{spm} * \\ 0 \end{spm}, \begin{spm} * & * \end{spm} \right)\\ \hline
(234) & 2 &  \left( \begin{spm} * \\ * \end{spm}, \begin{spm} 0 & * \end{spm} \right)\\ \hline
\begin{matrix} (123) \\ (13) \\ (1234) \\ (134) \end{matrix} & \begin{matrix} 2 \\ 3 \\ 3 \\ 4 \end{matrix} & \left( \begin{spm} 0 \\ 0 \end{spm}, \begin{spm} * & * \end{spm} \right)\\ \hline
\begin{matrix} (243) \\ (24) \\ (1432) \\ (142) \end{matrix} & \begin{matrix} 2 \\ 3 \\ 3 \\ 4 \end{matrix} & \left( \begin{spm} * \\ * \end{spm}, \begin{spm} 0 & 0 \end{spm} \right)\\ \hline
\begin{matrix} (12)(34) \\ (1243) \end{matrix} & \begin{matrix} 2 \\ 3  \end{matrix} & \left( \begin{spm} 0 \\ * \end{spm}, \begin{spm} * & 0 \end{spm} \right)\\ \hline
\begin{matrix} (1342) \\ (13)(24)\end{matrix} &  \begin{matrix}  3 \\ 4 \end{matrix} &  \left( \begin{spm} * \\ 0 \end{spm}, \begin{spm} 0 & * \end{spm} \right)\\ \hline
(1324) & 5 & \left( \begin{spm} 0 \\ 0 \end{spm}, \begin{spm} 0 & * \end{spm} \right)\\ \hline
(124) & 4 & \left( \begin{spm} 0 \\ 0 \end{spm}, \begin{spm} * & 0 \end{spm} \right) \\ \hline
(143) & 4 & \left( \begin{spm} 0 \\ * \end{spm}, \begin{spm} 0 & 0 \end{spm} \right) \\ \hline
(1423) &5  & \left( \begin{spm} * \\ 0 \end{spm}, \begin{spm} 0 & 0 \end{spm} \right)\\ \hline
\begin{matrix} (14) \\ (14)(23) \end{matrix} & \begin{matrix} 5 \\ 6 \end{matrix}  &  \left( \begin{spm} 0 \\ 0 \end{spm}, \begin{spm} 0 & 0 \end{spm} \right)
\end{array}
\end{equation*}
\end{table}
\mbox{ } \\
Using this table and the explicit description of the $R_{\textbf{i}, \textbf{a}_{\textbf{c}}}$ in \eqref{eq: spaces R i a for v eq 121}, we can compute the intersections
\begin{equation*}
R_{\textbf{i}, \textbf{a}_{\textbf{c}}}  \cap {}^w \mathfrak{n} =  R_{\textbf{i}, \textbf{a}_{\textbf{c}}}  \cap (  R_{\textbf{v}} \cap {}^w \mathfrak{n}).
\end{equation*}
We now have all the ingredients needed to evaluate the formula from \cref{lem: formula for h}. For example, for
\begin{align*}
\textbf{c}_4 & = (0,1,0,0,1,0) \\
\textbf{a}_{\textbf{c}_4} & = (0,1,1,1,1,0)\\
R_{\textbf{i}, \textbf{a}_{\textbf{c}_4}} &= \left( \begin{spm} 0 \\ * \end{spm} , \begin{spm} * & 0 \end{spm} \right)
\end{align*}
we get
\begin{align*}
h_{\textbf{c}_4}(v) & = v^{- 2} (1-v^{-2})^{-4} ( v^4 + v^6 + v^6 + v^6 + v^6 + v^6 + v^6 +v^8 +v^8+ v^{10} +  v^6 +v^8 +v^8\\
& \quad  + v^{10} + v^8+ v^{10} +v^6+ v^8 + v^{10} + v^{10}+ v^{10}+ v^{10} + v^{10} + v^{12}) \\
& = \tfrac{ v^2 + 8v^4 + 6v^6 + 8v^{8} + v^{10}}{(1-v^{-2})^4 }.
\end{align*}
More generally, we get
\begin{equation*}
 H = (1-v^{-2})^{-4} \begin{spm} v+2v^3+3v^5+3v^7+2v^9+v^{11}
 \\ v+5v^3+6v^5+6v^7+5v^9+v^{11}
 \\ 3v^3+3v^5+3v^7+3v^9
 \\ v^2+8v^4+6v^6+8v^8+v^{10}
\\ 3v^3+9v^5+9v^7+3v^9
\end{spm}.
\end{equation*}
We have already computed $Q$ and $\Psi$ in \eqref{eq: psi in the case 121} and \eqref{eq: matrix p and q in the v is 121 example} and hence we can compute
\begin{equation*}
F = Q^T \Psi^{-1} H = \begin{spm} 3v^9 +v^{11} \\  2v^5+2v^7+2v^9 \\ 2v^5+2v^7+2v^9 \\ v^6+v^8 \\ v^3 + 3v^5 \end{spm}.
\end{equation*}
Setting $v = 1$, we get that the dimensions of the simple $\mathcal{H}^{\aff}$-modules with central character $\chi_{(s,q_0)}$ are
\begin{equation*}
\dim L_1 = 4,\quad \dim L_2 = 6, \quad \dim L_3 = 6, \quad \dim L_4 = 2, \quad \dim L_5 = 4.
\end{equation*}
\end{example}

\bibliographystyle{alpha}
\bibliography{bibliography}

\end{document}